\documentclass[11pt]{article}
\usepackage{amsmath, upgreek}
\usepackage{amssymb}
\usepackage{amscd}
\usepackage{xspace}
\usepackage{verbatim}

\newcommand{\mysection}[1]{
\section{#1}\setcounter{equation}{0}}
\title{\bf Nonlinear boundary value problems \\relative to harmonic functions}
\author{{\bf Y.Oussama Boukarabila\footnote{Laboratoire d\'{}Analyse Nonlin´eaire et Math\'ematiques Appliqu\'{e}es, D\'{e}partement de Math\'ematiques,
			Universit\'e de Tlemcen, Algerie. Email: boukarabila.youcef.oussama@gmail.com}}\\[2mm]
 {\bf Laurent V\'eron\footnote{Institut Denis Poisson, Universit\'e de Tours, France. Email: veronl@univ-tours.fr}}\\[2mm]
}

\date{}
\begin{document}
 \maketitle


\newcommand{\txt}[1]{\;\text{ #1 }\;}
\newcommand{\tbf}{\textbf}
\newcommand{\tit}{\textit}
\newcommand{\tsc}{\textsc}
\newcommand{\trm}{\textrm}
\newcommand{\mbf}{\mathbf}
\newcommand{\mrm}{\mathrm}
\newcommand{\bsym}{\boldsymbol}
\newcommand{\scs}{\scriptstyle}
\newcommand{\sss}{\scriptscriptstyle}
\newcommand{\txts}{\textstyle}
\newcommand{\dsps}{\displaystyle}
\newcommand{\fnz}{\footnotesize}
\newcommand{\scz}{\scriptsize}
\newcommand{\be}{\begin{equation}}
\newcommand{\bel}[1]{\begin{equation}\label{#1}}
\newcommand{\ee}{\end{equation}}
\newcommand{\eqnl}[2]{\begin{equation}\label{#1}{#2}\end{equation}}
\newcommand{\barr}{\begin{eqnarray}}
\newcommand{\earr}{\end{eqnarray}}
\newcommand{\bars}{\begin{eqnarray*}}
\newcommand{\ears}{\end{eqnarray*}}
\newcommand{\nnu}{\nonumber \\}
\newtheorem{subn}{\name}
\renewcommand{\thesubn}{}
\newcommand{\bsn}[1]{\def\name{#1}\begin{subn}}
\newcommand{\esn}{\end{subn}}
\newtheorem{sub}{\name}[section]
\newcommand{\dn}[1]{\def\name{#1}}   
\newcommand{\bs}{\begin{sub}}
\newcommand{\es}{\end{sub}}
\newcommand{\bsl}[1]{\begin{sub}\label{#1}}
\newcommand{\bth}[1]{\def\name{Theorem}
\begin{sub}\label{t:#1}}
\newcommand{\blemma}[1]{\def\name{Lemma}
\begin{sub}\label{l:#1}}
\newcommand{\bcor}[1]{\def\name{Corollary}
\begin{sub}\label{c:#1}}
\newcommand{\bdef}[1]{\def\name{Definition}
\begin{sub}\label{d:#1}}
\newcommand{\bprop}[1]{\def\name{Proposition}
\begin{sub}\label{p:#1}}
\newcommand{\R}{\eqref}
\newcommand{\rth}[1]{Theorem~\ref{t:#1}}
\newcommand{\rlemma}[1]{Lemma~\ref{l:#1}}
\newcommand{\rcor}[1]{Corollary~\ref{c:#1}}
\newcommand{\rdef}[1]{Definition~\ref{d:#1}}
\newcommand{\rprop}[1]{Proposition~\ref{p:#1}}
\newcommand{\BA}{\begin{array}}
\newcommand{\EA}{\end{array}}
\newcommand{\BAN}{\renewcommand{\arraystretch}{1.2}
\setlength{\arraycolsep}{2pt}\begin{array}}
\newcommand{\BAV}[2]{\renewcommand{\arraystretch}{#1}
\setlength{\arraycolsep}{#2}\begin{array}}
\newcommand{\BSA}{\begin{subarray}}
\newcommand{\ESA}{\end{subarray}}
\newcommand{\BAL}{\begin{aligned}}
\newcommand{\EAL}{\end{aligned}}
\newcommand{\BALG}{\begin{alignat}}
\newcommand{\EALG}{\end{alignat}}
\newcommand{\BALGN}{\begin{alignat*}}
\newcommand{\EALGN}{\end{alignat*}}
\newcommand{\note}[1]{\textit{#1.}\hspace{2mm}}
\newcommand{\Proof}{\note{Proof}}
\newcommand{\qeda}{\hspace{10mm}\hfill $\square$}
\newcommand{\qed}{\\
${}$ \hfill $\square$}
\newcommand{\Remark}{\note{Remark}}
\newcommand{\modin}{$\,$\\[-4mm] \indent}
\newcommand{\forevery}{\quad \forall}
\newcommand{\set}[1]{\{#1\}}
\newcommand{\setdef}[2]{\{\,#1:\,#2\,\}}
\newcommand{\setm}[2]{\{\,#1\mid #2\,\}}
\newcommand{\mt}{\mapsto}
\newcommand{\lra}{\longrightarrow}
\newcommand{\lla}{\longleftarrow}
\newcommand{\llra}{\longleftrightarrow}
\newcommand{\Lra}{\Longrightarrow}
\newcommand{\Lla}{\Longleftarrow}
\newcommand{\Llra}{\Longleftrightarrow}
\newcommand{\warrow}{\rightharpoonup}
\newcommand{
\paran}[1]{\left (#1 \right )}
\newcommand{\sqbr}[1]{\left [#1 \right ]}
\newcommand{\curlybr}[1]{\left \{#1 \right \}}
\newcommand{\abs}[1]{\left |#1\right |}
\newcommand{\norm}[1]{\left \|#1\right \|}
\newcommand{
\paranb}[1]{\big (#1 \big )}
\newcommand{\lsqbrb}[1]{\big [#1 \big ]}
\newcommand{\lcurlybrb}[1]{\big \{#1 \big \}}
\newcommand{\absb}[1]{\big |#1\big |}
\newcommand{\normb}[1]{\big \|#1\big \|}
\newcommand{
\paranB}[1]{\Big (#1 \Big )}
\newcommand{\absB}[1]{\Big |#1\Big |}
\newcommand{\normB}[1]{\Big \|#1\Big \|}
\newcommand{\produal}[1]{\langle #1 \rangle}

\newcommand{\thkl}{\rule[-.5mm]{.3mm}{3mm}}
\newcommand{\thknorm}[1]{\thkl #1 \thkl\,}
\newcommand{\trinorm}[1]{|\!|\!| #1 |\!|\!|\,}
\newcommand{\bang}[1]{\langle #1 \rangle}
\def\angb<#1>{\langle #1 \rangle}
\newcommand{\vstrut}[1]{\rule{0mm}{#1}}
\newcommand{\rec}[1]{\frac{1}{#1}}
\newcommand{\opname}[1]{\mbox{\rm #1}\,}
\newcommand{\supp}{\opname{supp}}
\newcommand{\dist}{\opname{dist}}
\newcommand{\myfrac}[2]{{\displaystyle \frac{#1}{#2} }}
\newcommand{\myint}[2]{{\displaystyle \int_{#1}^{#2}}}
\newcommand{\mysum}[2]{{\displaystyle \sum_{#1}^{#2}}}
\newcommand {\dint}{{\displaystyle \myint\!\!\myint}}
\newcommand{\q}{\quad}
\newcommand{\qq}{\qquad}
\newcommand{\hsp}[1]{\hspace{#1mm}}
\newcommand{\vsp}[1]{\vspace{#1mm}}
\newcommand{\ity}{\infty}
\newcommand{\prt}{\partial}
\newcommand{\sms}{\setminus}
\newcommand{\ems}{\emptyset}
\newcommand{\ti}{\times}
\newcommand{\pr}{^\prime}
\newcommand{\ppr}{^{\prime\prime}}
\newcommand{\tl}{\tilde}
\newcommand{\sbs}{\subset}
\newcommand{\sbeq}{\subseteq}
\newcommand{\nind}{\noindent}
\newcommand{\ind}{\indent}
\newcommand{\ovl}{\overline}
\newcommand{\unl}{\underline}
\newcommand{\nin}{\not\in}
\newcommand{\pfrac}[2]{\genfrac{(}{)}{}{}{#1}{#2}}

\def\ga{\alpha}     \def\gb{\beta}       \def\gg{\gamma}
\def\gc{\chi}       \def\gd{\delta}      \def\ge{\epsilon}
\def\gth{\theta}                         \def\vge{\varepsilon}
\def\gf{\phi}       \def\vgf{\varphi}    \def\gh{\eta}
\def\gi{\iota}      \def\gk{\kappa}      \def\gl{\lambda}
\def\gm{\mu}        \def\gn{\nu}         \def\gp{\pi}
\def\vgp{\varpi}    \def\gr{\rho}        \def\vgr{\varrho}
\def\gs{\sigma}     \def\vgs{\varsigma}  \def\gt{\tau}
\def\gu{\upsilon}   \def\gv{\vartheta}   \def\gw{\omega}
\def\gx{\xi}        \def\gy{\psi}        \def\gz{\zeta}
\def\Gg{\Gamma}     \def\Gd{\Delta}      \def\Gf{\Phi}
\def\Gth{\Theta}
\def\Gl{\Lambda}    \def\Gs{\Sigma}      \def\Gp{\Pi}
\def\Gw{\Omega}     \def\Gx{\Xi}         \def\Gy{\Psi}

\def\CS{{\mathcal S}}   \def\CM{{\mathcal M}}   \def\CN{{\mathcal N}}
\def\CR{{\mathcal R}}   \def\CO{{\mathcal O}}   \def\CP{{\mathcal P}}
\def\CA{{\mathcal A}}   \def\CB{{\mathcal B}}   \def\CC{{\mathcal C}}
\def\CD{{\mathcal D}}   \def\CE{{\mathcal E}}   \def\CF{{\mathcal F}}
\def\CG{{\mathcal G}}   \def\CH{{\mathcal H}}   \def\CI{{\mathcal I}}
\def\CJ{{\mathcal J}}   \def\CK{{\mathcal K}}   \def\CL{{\mathcal L}}
\def\CT{{\mathcal T}}   \def\CU{{\mathcal U}}   \def\CV{{\mathcal V}}
\def\CZ{{\mathcal Z}}   \def\CX{{\mathcal X}}   \def\CY{{\mathcal Y}}
\def\CW{{\mathcal W}} \def\CQ{{\mathcal Q}}
\def\BBA {\mathbb A}   \def\BBb {\mathbb B}    \def\BBC {\mathbb C}
\def\BBD {\mathbb D}   \def\BBE {\mathbb E}    \def\BBF {\mathbb F}
\def\BBG {\mathbb G}   \def\BBH {\mathbb H}    \def\BBI {\mathbb I}
\def\BBJ {\mathbb J}   \def\BBK {\mathbb K}    \def\BBL {\mathbb L}
\def\BBM {\mathbb M}   \def\BBN {\mathbb N}    \def\BBO {\mathbb O}
\def\BBP {\mathbb P}   \def\BBR {\mathbb R}    \def\BBS {\mathbb S}
\def\BBT {\mathbb T}   \def\BBU {\mathbb U}    \def\BBV {\mathbb V}
\def\BBW {\mathbb W}   \def\BBX {\mathbb X}    \def\BBY {\mathbb Y}
\def\BBZ {\mathbb Z}   \def\BBQ {\mathbb Q}

\def\GTA {\mathfrak A}   \def\GTB {\mathfrak B}    \def\GTC {\mathfrak C}
\def\GTD {\mathfrak D}   \def\GTE {\mathfrak E}    \def\GTF {\mathfrak F}
\def\GTG {\mathfrak G}   \def\GTH {\mathfrak H}    \def\GTI {\mathfrak I}
\def\GTJ {\mathfrak J}   \def\GTK {\mathfrak K}    \def\GTL {\mathfrak L}
\def\GTM {\mathfrak M}   \def\GTN {\mathfrak N}    \def\GTO {\mathfrak O}
\def\GTP {\mathfrak P}   \def\GTR {\mathfrak R}    \def\GTS {\mathfrak S}
\def\GTT {\mathfrak T}   \def\GTU {\mathfrak U}    \def\GTV {\mathfrak V}
\def\GTW {\mathfrak W}   \def\GTX {\mathfrak X}    \def\GTY {\mathfrak Y}
\def\GTZ {\mathfrak Z}   \def\GTQ {\mathfrak Q}

\font\Sym= msam10 
\def\SYM#1{\hbox{\Sym #1}}
\newcommand{\bdw}{\prt\Gw\xspace}

\date{}
\maketitle
\begin {center}
{\it
 To Shair with high esteem and sincere friendship}
\end{center}\medskip

\abstract { We study the problem of finding a function $u$ verifying $-\Gd u=0$ in $\Gw $ under the boundary condition $\frac{\prt u}{\prt {\bf n}}+g(u)=\gm$ on $\prt\Gw$ where $\Gw\subset\BBR^N$  is a smooth domain, ${\bf n}$ the normal unit outward vector to $\Gw$, $\gm$ is a measure on $\prt\Gw$ and $g$ a continuous nondecreasing function. We give sufficient condition on $g$ for this problem to be solvable for any measure. When 
$g(r)=|r|^{p-1}r$, $p>1$, we give conditions in order an isolated singularity on $\prt\Gw$ be removable. We also give capacitary conditions on a measure $\gm$ in order the problem with $g(r)=|r|^{p-1}r$ to be solvable for some $\gm$. We also study the isolated singularities of functions satisfying $-\Gd u=0\quad\text{in }\Gw$ and  $\frac{\prt u}{\prt {\bf n}}+g(u)=0$ on $\prt\Gw\setminus\{0\}$.\smallskip

  \vspace{1mm}

\tableofcontents
\vspace{1mm}
\hspace{.05in}

\tableofcontents
\mysection{Introduction}
Let $\Gw$ be a smooth bounded domain in $\BBR^N$ such that $0\in\prt\Gw$ and $g:\BBR\mapsto\BBR$ a continuous nondecreasing  function such that 
$rg(r)\geq 0$. The aim of this article is to study the following nonlinear problem
\bel{I-1-0}\BA {lll}
\phantom{g(u;}
-\Gd u+u=0\qquad&\text{in }\;\Gw\\[0mm]\phantom{,,,u}
\myfrac{\prt u}{\prt {\bf n}}+g(u)=\gm\qquad&\text{in }\;\prt\Gw,
\EA
\ee
where $\gm$ is a Radon measure on $\prt\Gw$ and ${\bf n}$ the outward normal unit vector on $\prt\Gw$. An associated model problem on which we can develop sharp estimate is the following equation in the upper half-space $\BBR^N_+:=\{x=(x_1,...,x_{_N})\in\BBR^N:x_N>0\}$,
\bel{I-1-2}\BA {lll}
\phantom{g(u-,,,}
-\Gd u=0\qquad&\text{in }\;\BBR^N_+\\[0mm]
-\myfrac{\prt u}{\prt x_{_N}}\!+|u|^{p-1}u=0\qquad&\text{in }\;\prt\BBR^N_+\setminus\{0\}
\EA
\ee
 where $p>1$. These two problems are by essence non-local and actually, the second problem can be expressed by introducing the square root of the Laplacian in $\BBR^{N-1}$ under the form
 \bel{I-1-3}\BA {lll}
\phantom{g(u---}
(-\Gd_{_{N-1}})^{\frac{1}{2}}\tilde u+\tilde u^p=0\qquad\text{in }\;\BBR^{N-1}\setminus\{0\},
\EA
\ee
where $\displaystyle\Gd_{_{N-1}}=\sum_{j=1}^{N-1}\frac{\prt^2}{\prt x_j^2}$ and $\tilde u(x_1,...,x_{N-1})=u(x_1,...,x_{N-1},0)$. The second equation is equivariant under the scaling transformation $T_k$ ($k>0$) defined by
 \bel{I-1-4}\BA {lll}
T_k[u](x)=k^{\frac{1}{p-1}}u(kx).
\EA
\ee
Therefore it is natural to look for self-similar solutions i.e. solutions satisfying $T_k[u]=u$ for any $k>0$. Introducing the spherical coordinates $(r,\gs)\in (0,\infty\ti S^{N-1})$, then a self-similar solution endows the form
 \bel{I-1-5}\BA {lll}
u(x)=u(r,\gs)=r^{-\frac{1}{p-1}}\gw(\gs),
\EA
\ee
and $\gw$ satisfies
 \bel{I-1-6}\BA {lll}
\!\phantom{,,}\Gd'\gw+\ell_{N,p}\gw=0\qquad&\text{in }\;S^{N-1}_+\\
\myfrac{\prt \gw}{\prt \gn}+|\gw|^{p-1}\gw=0\qquad&\text{in }\;\prt S^{N-1}_+,
\EA
\ee
where $\Gd'$ is the Laplace-Beltrami operator on the unit sphere $S^{N-1}$, $\gn$ is the outward normal unit vector to $\prt S^{N-1}_+$ tangent to $S^{N-1}$ and 
 \bel{I-1-7}\BA {lll}
\ell_{N,p}=\left(\myfrac{1}{p-1}\right)\left(\myfrac{1}{p-1}+2-N\right). 
\EA
\ee
This problem points out the existence of critical values of $p$. We denote by $\CE$ the set of solutions of $(\ref{I-1-6})$ and $\CE_+=\{\gw\in\CE:\gw\geq 0\}$. This set has the following structure: \medskip

\noindent{\bf Theorem A}. {\it 1- If
 \bel{I-1-8c}\BA {lll}
 p\geq\myfrac{N-1}{N-2},
\EA
\ee
then $\CE=\{0\}$. \smallskip

\noindent  2- If 
 \bel{I-1-8a}\BA {lll}
 1<p\leq \myfrac{N}{N-1},
\EA
\ee
then $\CE_+=\{0\}$. 
\smallskip

\noindent 3- If
 \bel{I-1-8b}\BA {lll}
\myfrac{N}{N-1}< p<\myfrac{N-1}{N-2},
\EA
\ee
then $\CE=\{\gw_s,-\gw_s,0\}$ where $\gw_s$ is the unique positive solution of $(\ref{I-1-6})$.}\medskip

\noindent When $1<p<\frac{N}{N-1}$, we show that there exist signed solutions to $(\ref{I-1-6})$. \medskip 


\noindent{\bf Theorem B}. {\it Let $\Gw\subset\BBR^N$, $N\geq 2$, be a bounded $C^2$ domain such that $0\in\prt\Gw$ and 
 $g:\BBR\mapsto\BBR$ a continuous function which satisfies $sg(s)\geq 0$. Then any function 
 $u\in C^1(\overline\Gw\setminus\{0\})$ solution of 
 \bel{I-1-10}\BA {lll}
\phantom{g(,,}
-\Gd u=0\qquad&\text{in }\;\Gw\\[0mm]
\!\myfrac{\prt u}{\prt {\bf n}}+g(u)=0\qquad&\text{in }\;\partial\Gw\setminus\{0\},
\EA
\ee
satisfying near $x=0$ either $u(x)=o(|x|^{2-N})$  if  $N\geq 3$ or $u(x)=o(\ln|x|)$ if $N=2$, is constant and $g(u)=0$.
  }\medskip

The set $\CE$ plays a fundamental role in the characterization of boundary isolated singularities of solutions of 
\bel{I-1-9}\BA {lll}
\phantom{g(u-,}
-\Gd u=0\qquad&\text{in }\;\Gw\\[0mm]
\!\!\!\myfrac{\prt u}{\prt {\bf n}}+|u|^{p-1}u=0\qquad&\text{in }\;\partial\Gw\setminus\{0\}.
\EA
\ee

 \noindent{\bf Theorem C}. {\it Let $\Gw\subset\BBR^N$ be a smooth bounded domain such that $0\in\prt\Gw$. Assume $u\in C^1(\overline\Gw\setminus\{0\})$ is a nonnegative function satisfying $(\ref{I-1-9})$ and such that $|x|^{\frac{1}{p-1}}u(x)$ is bounded.\smallskip
 
 \noindent 1- If $p\geq\frac{N-1}{N-2}$, then $u=0$\smallskip
 
 \noindent 2- If $\frac{N}{N-1}<p<\frac{N-1}{N-2}$, we have the following alternative:
 
 \noindent 2-(i)  either 
    \bel{I-1-14}\BA {lll}\displaystyle
 \lim_{r\to 0} r^{\frac{1}{p-1}}u(r,\gs)=\gw_s(\gs)\quad\text{locally uniformly on } S^{N-1}_+,
 \EA \ee
\smallskip
 
  \noindent 2-(ii) or there exists a nonnegative real number $k$ such that there holds,
    \bel{I-1-13}\BA {lll}\displaystyle
a)\qquad\qquad \qquad&\displaystyle\lim_{|x|\to 0} |x|^{2-N}u(x)=k\quad&\text{if }\; N\geq 3,\qquad\qquad\\[1mm]
\displaystyle 
b)\displaystyle\qquad\qquad &\displaystyle\lim_{|x|\to 0} (-\ln |x|)^{-1}u(x)=k\quad&\text{if }\; N=2.\qquad\qquad
\EA \ee
\smallskip
 }\medskip
 
  \noindent   The assumption on the boundedness of $|x|^{\frac{1}{p-1}}u(x)$ seems necessary since no Keller-Osserman universal estimate \cite{Ke}, \cite{Os} appears to hold. Actually, if 
 $u$ satisfies $(\ref{I-1-2})$, the function $\tilde u$ defined in whole $\BBR^N$ by 
 \bel{I-1-9'}
 \tilde u(x_1,...,x_{_N})=\left\{\BA {lll}u(x_1,...,x_{_N})\quad&\text{if }x_{_N}>0\\
 u(x_1,...,-x_{_N})\quad&\text{if }x_{_N}<0,
\EA \right.
 \ee
 satisfies
  \bel{I-1-9''}
-\Gd \tilde u+2|u|^{p-1}u\CH_{\prt\BBR^N_+}=0\qquad\text{in }\BBR^N_+\setminus\{0\},
 \ee
 where $\CH_{\prt\BBR^N_+}$ is the (N-1)-dimensional Lebesgue measure supported by $\prt\BBR^N_+$. Hence the coercivity due to the nonlinear term is 
 localized on $\prt\BBR^N_+$. Such problems with measure valued nonlinear potential are studied in \cite{SainVer}. Notice also that when $p>\frac{N-1}{N-2}$, then the assumption $u(x)=O(|x|^{-\frac{1}{p-1}})$ implies that $u(x)=o(|x|^{2-N})$, hence Theorem B implies Theorem C.
\medskip

When $u$ satisfies $(\ref{I-1-13})$, the problem can be  interpreted with a boundary data holding in the sense of distributions,
 \bel{I-1-15}\BA {lll}
\phantom{u^p}
-\Gd u=0\qquad&\text{in }\;\Gw\\[0mm]
\myfrac{\prt u}{\prt {\bf n}}+u^p=k\gd_0\qquad&\text{in }\;\CD'(\prt\Gw).
\EA
\ee
 For more general measures and nonlinearities, we define a solution of problem $(\ref{I-1-0})$ as follows,
  \medskip
 
 \noindent{\bf Definition}. {\it Let $\Gw\subset\BBR^N$ be as in Theorem B , $\gm\in\mathfrak M(\prt\Gw)$ and $g:\BBR\mapsto\BBR$ be a continuous function. A function $u\in L^1(\Gw)$  is a weak solution of $(\ref{I-1-0})$ if it admits a boundary trace $u\lfloor_{\prt\Gw}$ which is a Borel function on $\prt\Gw$, $g(u)\in L^1(\prt\Gw)$ and 
  \bel{I-1-16}\BA {lll}
\myint{\Gw}{}u(-\Gd\xi+\xi) dx+\myint{\prt\Gw}{}g(u)\xi dS=\myint{\prt\Gw}{}\xi d\gm,\quad\text{for all }\;\xi\in \CC(\Gw),
\EA
\ee
where 
 \bel{I-1-17}\CC(\Gw)=\left\{\xi\in C^1(\overline\Gw):\Gd\xi\in L^{\infty}(\Gw),\,\myfrac{\prt \xi}{\prt {\bf n}}=0\;\text{ on }\prt\Gw\right\}.\ee
}\medskip

In the next result we give a condition for the unconditionnal solvability of problem $(\ref{I-1-0})$. \medskip
 
 \noindent{\bf Theorem D}. {\it Let $\Gw\subset\BBR^N$, $N\geq 3$ be a bounded $C^2$ domain and $g:\BBR\mapsto\BBR$ a continuous nondecreasing function such that $g(0)=0$. If $g$ satisfies 
  \bel{I-1-20}\BA {lll}
\myint{1}{\infty}(g(s)+|g(-s)|)s^{-\frac{2N-3}{N-2}}ds<\infty,
\EA
\ee
then for any $\gm\in\mathfrak M(\prt\Gw)$, the problem $(\ref{I-1-0})$ admits a unique solution.
}\medskip

A nonlinearity which satisfies $(\ref{I-1-20})$ is called {\it subcritical}. When $N=2$ this notion has to be modified. Following V\`azquez we define the exponential orders of growth of a continuous nondecreasing function $g:\BBR\mapsto\BBR$ vanishing at $0$ by 
  \bel{I-1-21}a_+(g)=\inf\left\{a\geq 0:\myint{0}{\infty}e^{-as}g(s) ds<\infty\right\},
\ee
and 
  \bel{I-1-22}a_-(g)=sup\left\{a\leq 0:\myint{-\infty}{0}e^{as}g(s) ds>-\infty\right\}.
\ee
 \medskip
 
 \noindent{\bf Theorem E}. {\it Let $\Gw\subset\BBR^2$ be a bounded $C^2$ domain and $g:\BBR\mapsto\BBR$ a continuous nondecreasing function such that $g(0)=0$. \smallskip
 
\noindent 1- If $a_+(g)=a_-(g)=0$, then for any $\gm\in\mathfrak M(\prt\Gw)$ the problem $(\ref{I-1-0})$ admits a unique solution,
\smallskip
 
\noindent 2- if $0<a_+(g)<\infty $ and $-\infty<a_-(g)<0$ the problem $(\ref{I-1-0})$  admits a unique solution with $\gm=\displaystyle\sum_{j=1}^k\ga_j\gd_{a_j}$, with 
$a_j\in\prt\Gw$ and $\ga_j\in\BBR^*$, provided
\bel{-1-22}
\myfrac{\gp}{a_-(g)}\leq \ga_j\leq \myfrac{\gp}{a_+(g)}.
\ee
}\medskip

When $N\geq 3$ and $g$ does not satisfy $(\ref{I-1-20})$, there may not exist solutions for any measure. The problem is well understood if $g(r)=|r|^{p-1}r$. For example, if $p\geq \frac{N-1}{N-2}$, there is no weak solution to the problem
  \bel{I-1-23}\BA {lll} 
\phantom{|u^{p-1}u,}
-\Gd u+u=0\qquad&\text{in }\;\Gw\\[0mm]\phantom{+,u}
\myfrac{\prt u}{\prt {\bf n}}+ |u|^{p-1}u=\gm\qquad&\text{in }\;\CD'(\prt\Gw).
\EA
\ee
when $\gm=\ga\gd_a$ with $a\in\prt\Gw$. As in many similar problems, the condition for a Radon measure in order there exists a weak solution to $(\ref{I-1-23})$ is expressed in terms of Bessel capacities, presently the capacity $C^{1,p'}_{{\prt\Gw}}$ on the boundary with $p'=\frac{p}{p-1}$.\medskip

 \noindent{\bf Theorem F}. {\it Let $\Gw\subset\BBR^N$, $N\geq 2$ be a bounded $C^2$. Then problem $(\ref{I-1-23})$ admits a solution with $\gm\in\mathfrak M_+(\prt\Gw)$, necessarily unique, if and only if $\gm$ vanishes on Borel set $E\subset\prt\Gw$ such that $C^{1,p'}_{{\prt\Gw}}(E)=0$.
 }\medskip
 
 This work is the main part of the PhD thesis of the first author prepared in the {\it Laboratoire de Math\'ematiques et Physique Th\'eorique} of the University of Tours under the supervision of the second author.

\mysection{Separable solutions}
We recall that the upper hemisphere $S^{N-1}_+$ can be parametrized as follows
 \bel{I-2-1}S^{N-1}_+=\{\gs=(\sin\gf\, \gs,\cos\gf):\gs'\in S^{N-2}, \gf\in [0,\tfrac\gp 2]\},
\ee
and we write $\gw(\gs)=\gw(\gs',\gf))$. With this parametrization the Laplace-Beltrami operator on $S^{N-1}$ endows the form
 \bel{I-2-2}\Gd'\gw=\myfrac{1}{\sin^{N-2}\gf}\left(\sin^{N-2}\gf\,\gw_\gf\right)_\gf+\myfrac{1}{\sin^2\gf}\Gd''\gw
\ee
where $\Gd''$ is the Laplace-Beltrami operator on $S^{N-2}$. The surface measure on $S^{N-1}$ induced by the Euclidean metric in $\BBR^N$ is $dS(\gs)=\sin^{N-2}\gf dS'(\gs')d\gf$ where $dS'(\gs')$ is the surface measure on $S^{N-2}$ induced by the Euclidean metric in $\BBR^{N-1}$. 

\subsection{Proof of Theorem {\bf A}}
\noindent {\it Proof of assertion 1.} If $p\geq \frac{N-1}{N-2}$ then $\ell_{N,p}\leq 0$. If $\gw$ is a solution of $(\ref{I-1-6})$, then
$$\myint{S^{N-1}_+}{}\left(|\nabla'\gw|^2-\ell_{N,p}\gw^2\right)dS+\myint{\prt S^{N-1}_+}{}|\gw|^{p+1}dS'=0.
$$ 
Hence $\gw=0$. \smallskip

\noindent {\it Proof of assertion  2.} Assume $(\ref{I-1-8a})$ holds and $\gw$ is a positive solution of $(\ref{I-1-6})$. The function $\phi\mapsto\cos\gf$ is the first 
eigenfunction of $-\Gd'$ in $H^1_0(S^{N-1}_+)$ with corresponding eigenvalue N-1. Multiplying the equation by $\cos\gf$ and integrating yields
$$\left(\ell_{N,p}+1-N\right)\myint{S^{N-1}_+}{}\gw\cos\gf dS+\myint{\prt S^{N-1}_+}{}\gw dS'=0
$$
If $1<p\leq \frac{N}{N-1}$, then $\ell_{N,p}+1-N\geq 0$, hence $\gw\lfloor_{\prt S^{N-1}_+}=0$, hence $\gw=0$ by Hopf boundary lemma..\smallskip

\noindent {\it Proof of assertion  3.} Assume $(\ref{I-1-8b})$ holds. We first prove that any solution $\gw$ of $(\ref{I-1-6})$ depends only on $\gf$  following a method introduced in \cite{Ver} and it has constant sign.
We set
$$\bar\gw(\gf)=\myfrac{1}{|S^{N-2}|}\myint{S^{N-2}}{}\gw(\gs',\gf)dS'(\gs').
$$
Then
$$\BA{lll}
\myint{S^{N-2}}{}\left(|\gw|^{p-1}\gw-\overline{|\gw|^{p-1}\gw}\right)\left(\gw-\overline{\gw}\right)dS'=
\myint{S^{N-2}}{}\left(|\gw|^{p-1}\gw-{|\overline\gw|^{p-1}\overline\gw}\right)\left(\gw-\overline{\gw}\right)dS',
\EA$$
since
$$\myint{S^{N-2}}{}\!\!\left(\overline{|\gw|^{p-1}\gw}-{|\overline\gw|^{p-1}\overline\gw}\right)\left(\gw-\overline{\gw}\right)dS'=
\left(\overline{|\gw|^{p-1}\gw}-{|\overline\gw|^{p-1}\overline\gw}\right)\!\myint{S^{N-2}}{}\!\!\left(\gw-\overline{\gw}\right)dS'=0.
$$
Hence
$$\myint{S^{N-2}}{}\left(|\gw|^{p-1}\gw-\overline{|\gw|^{p-1}\gw}\right)\left(\gw-\overline{\gw}\right)dS'
\geq 2^{1-p}\myint{S^{N-2}}{}\left|\gw-\overline{\gw}\right|^{p+1}dS'.
$$
From the expression $(\ref{I-2-2})$ we get
$$\BA{lll}
-\myint{S^{N-1}_+}{}|\nabla'(\gw-\overline\gw)|^2dS+\ell_{N,p}
\myint{S^{N-1}_+}{}(\gw-\overline\gw)^2dS\\[4mm]
\phantom{---------}
=\myint{S^{N-2}}{}\left(|\gw|^{p-1}\gw-\overline{|\gw|^{p-1}\gw}\right)\left(\gw-\overline{\gw}\right)dS'
\\[4mm]
\phantom{---------}
\geq 2^{1-p}\myint{S^{N-2}}{}\left|\gw-\overline{\gw}\right|^{p+1}dS'.
\EA$$
Since $\overline\gw$ is the projection of $\gw$ onto the first eigenspace of $-\Gd'$ in $H^1(S^{N-1}_+)$ and N-1  the corresponding eigenvalue, 
$$-\myint{S^{N-1}_+}{}|\nabla'(\gw-\overline\gw)|^2dS\leq (1-N)\myint{S^{N-1}_+}{}(\gw-\overline\gw)^2dS.
$$
Hence
$$\left(\ell_{N,p}+1-N\right)\myint{S^{N-1}_+}{}(\gw-\overline\gw)^2dS\geq 2^{1-p}\myint{S^{N-2}}{}\left|\gw-\overline{\gw}\right|^{p+1}dS'.
$$
If $p\geq\frac{N}{N-1}$, then $\ell_{N,p}+1-N\leq 0$. This implies $\gw=\overline\gw$. It follows that $\gw$ depends only on the variable $\gf\in (0,\frac\gp2)$ and thus it satisfies
\bel{1-2-4}\BA{lll}
\myfrac{1}{\sin^{N-2}\gf}\left(\sin^{N-2}\gf\,\gw_\gf\right)_\gf+\ell_{N,p}\gw=0\qquad\text{in }\, (0,\frac{\gp}{2})\\
\gw_\gf(0)=0\,,\; \left(\gw_\gf+|\gw|^{p-1}\gw\right)\left(\frac\gp2\right)=0.
\EA\ee
Next we prove that any solution has constant sign. Let us assume that $\gw(0)>0$. If $\gw$ vanishes at a first point some $\gf_0\in (0,\frac{\gp}{2}]$, then it is positive on $(0,\gf_0)$ and $\gw_\gf(\gf_0)<0$ by Cauchy-Lipschitz theorem. If $\gf_0=\frac\gp2$, then $\gw_\gf(\gf_0)=0$ from $(\ref{1-2-4})$, contradiction. Hence $\gf_0<\frac\gp2$. This implies that 
$\gw$ is a positive solution of
\bel{1-2-5}\BA{lll}
\myfrac{1}{\sin^{N-2}\gf}\left(\sin^{N-2}\gf\,\gw_\gf\right)_\gf+\ell_{N,p}\gw=0\qquad\text{in }\, (0,\gf_0)\\
\gw_\gf(0)=0\,,\; \gw(\gf_0)=0.
\EA\ee
Thus $\gw$ is a first eigenfunction of $-\Gd'$ in $H^1_0(S_{\gf_0})$ where $$S_{\gf_0}=\{\gs=(\gs',\gf)\in S^{N-2}\ti (0,\gf_0)\}\subsetneq S^{N-1}_+.$$ Hence $\ell_{N,p}>N-1$, contradiction.\smallskip

\noindent Then we prove that there exists at most one positive solution $\gw$. Let $\tilde\gw$ be another positive solution. A straightformard computation yields
$$\BA {lll}
\displaystyle 0=\myint{S^{N-1}_+}{}\left(\myfrac{\Gd'\gw}{\gw}-\myfrac{\Gd'\tilde\gw}{\tilde\gw}\right)\left(\gw^2-\tilde\gw^2\right)dS\\[4mm]
\phantom{0}=
-\myint{S^{N-1}_+}{} \left(\myfrac{1}{\gw^2}+\myfrac{1}{\tilde\gw^2}\right)\left|\gw\nabla'\tilde\gw-\tilde\gw\nabla'\gw\right|^2
-\myint{\prt S^{N-1}_+}{} \left(\gw^{p-1}-\tilde\gw^{p-1}\right)\left(\gw^2-\tilde\gw^2\right) dS'.
\displaystyle 
\EA$$
This implies that $\gw=\tilde\gw$. \smallskip

\noindent Finally we prove existence. Set
 \bel{I-2-3}J(\eta)=\myfrac{1}{2}\myint{S^{N-1}_+}{}\left(|\nabla'\eta|^2-\ell_{N,p}\eta^2\right)dS+\myfrac{1}{p+1}\myint{\prt S^{N-1}_+}{}|\eta|^{p+1}dS.
 \ee
The functional  $J$ is defined in 
$$\BBX_{rad}(S^{N-1}_+):=\left\{\eta\in H^1(S^{N-1}_+)\cap L^{p+1}(\prt S^{N-1}_+):\eta\text { depends only on }\gf\in[0,\tfrac\gp2]\right\},$$
and it is lower continuous. If $\eta\in \BBX_{rad}(S^{N-1}_+)$, then
$\eta=\eta_1+\eta_0$ where $\eta_0=\eta(\frac\gp2)$ and  $\eta_1\in H_0^1(S^{N-1}_+)$. In particular
$$\BA {lll}J(\eta)=\myfrac{1}{2}\myint{S^{N-1}_+}{}\left(|\nabla'\eta_1|^2-\ell_{N,p}\eta_1^2\right)dS -\ell_{N,p}\eta_0\myint{S^{N-1}_+}{}\eta_1dS\\[4mm]
\phantom{---------------------}
-\myfrac{\ell_{N,p}|S^{N-1}_+|}{2}\eta_0^2
+\myfrac{|S^{N-2}|}{p+1}|\eta_0|^{p+1}
\EA$$
Since $(\ref{I-1-8b})$ holds, $0<\ell_{N,p}\leq N-1$; if we take $\eta(\gf)=\ge_0\in\BBR$, then 
$$\BA {lll}J(\eta)=-\myfrac{\ell_{N,p}|S^{N-1}_+|}{2}\ge_0^2
+\myfrac{|S^{N-2}|}{p+1}|\ge_0|^{p+1}.
\EA$$
Hence the infimum of $J$ in $\BBX_{rad}(S^{N-1}_+)$ is negative. Since $p>\frac{N}{N-1}$, then $\ell_{N,p}< N-1$, and
for $\ge= N-1-\ell_{N,p}>0$ there holds
$$J(\eta)\geq \myfrac{\ge}{2}\myint{S^{N-1}_+}{}\eta_1^2dS
-\ell_{N,p}\eta_0\myint{S^{N-1}_+}{}\eta_1dS\\[4mm]
-\myfrac{\ell_{N,p}|S^{N-1}_+|}{2}\eta_0^2
+\myfrac{|S^{N-2}|}{p+1}|\eta_0|^{p+1}.
$$
By Young's inequality $J(\eta)\to\infty$ when $\norm\eta_{H^1(S^{N-1}_+)}+\norm\eta_{L^{p+1}(\prt S^{N-1}_+)}\to\infty$. Therefore 
$J$ achieves its minimum in $\BBX_{rad}(S^{N-1}_+)$ at some $\gw$, which can be assume to be positive since $J(\eta)=J(|\eta|)$. If we denote it by $\gw_s$,
there holds $\CE=\{\gw_s,-\gw_s,0\}$, which ends the proof.\qeda

\medskip

The value $p=\frac{N}{N-1}$ is a bifurcation value as it is shown below.

\bprop{bif} There exists a $C^1$ curve $\ge\mapsto (p_\ge,\gw_\ge)$ defined in $[0,\ge_0]$ with $\ge>0$ such that $(p_0,\gw_0)=(\frac{N}{N-1},0)$ where $1<p_\ge<\frac{N}{N-1}$ and $\gw_\ge$ is a nonzero signed solution of 
 \bel{I-2-4}\BA {lll}
\!\phantom{,,}\!\!\Gd'\gw_\ge+\ell_{N,p_\ge}\gw=0\qquad&\text{in }\;S^{N-1}_+\\
\myfrac{\prt \gw}{\prt \gn}+|\gw|^{p_\ge-1}\gw=0\qquad&\text{in }\;\prt S^{N-1}_+.
\EA
\ee
\es
\Proof The linearization of $(\ref{I-1-6})$ at $p=\frac{N}{N-1}$ and $\gw=0$ yields 
 \bel{I-2-5}\BA {lll}
\Gd'\psi+(N-1)\psi=0\qquad&\text{in }\;S^{N-1}_+\\[1mm]\phantom{------}
\myfrac{\prt \psi}{\prt \gn}=0\qquad&\text{in }\;\prt S^{N-1}_+.
\EA
\ee
If $\BBR^N_+:=\{x=(x_1,...,x_{_N}):x_{_N}>0\}$, then for $j<N$ the restriction to $S^{N-1}_+$ of the function $\psi_j:x\mapsto x_j$ satisfies 
$(\ref{I-2-5})$. In order to satisfy the simplicity requirement, we consider the functions defined on $S^{N-1}_+$ depending only of the variable 
$x_j\lfloor_{S^{N-1}_+}$. Then $\psi_j$ is a simple eigenfunction of $\Gd'$ associated to the eigenvalue $N-1$. By the classical Crandall-Rabinowitz theorem\cite{CranRab} there exists a $C^1$ curve $\ge\mapsto (p_\ge,\gw_\ge)$ starting from $(\frac{N}{N-1},0)$ such that $\gw_\ge$ is a nonzero solution depending only of 
the variable $x_j\lfloor_{S^{N-1}_+}$ of the problem
 \bel{I-2-6}\BA {lll}
\phantom{-}\Gd'\gw_\ge+\ell_{N,p_\ge}\gw_\ge=0\qquad&\text{in }\;S^{N-1}_+\\[1mm]\phantom{}
\myfrac{\prt \gw_\ge}{\prt \gn}+|\gw_\ge|^{p_\ge-1}\gw_\ge=0\qquad&\text{in }\;\prt S^{N-1}_+.
\EA
\ee
Since $\gw_\ge$ depends only on $x_j\lfloor_{S^{N-1}_+}$ and inherits the properties of $\psi_j$, it changes sign. By Theorem A-1-2, $p_\ge<\frac{N}{N-1}$, which ends the proof.\qeda
\subsection{Separable solutions in dimension 2}
When $N=2$, $(\ref{1-2-5})$ endows the form
\bel{1-2-7}\BA{lll}
\phantom{--}\gw_{\gf\gf}+\myfrac{1}{(p-1)^2}\gw=0\qquad\text{in }\, (0,\gp)\\[3mm]
\left(-\gw_\gf+|\gw|^{p-1}\gw\right)\left(0\right)=0\\[1mm]
\phantom{0}\left(\gw_\gf+|\gw|^{p-1}\gw\right)\left(\gp\right)=0,
\EA\ee
and therefore
$$\gw(\gf)=a\cos\left(\frac{\gf}{p-1}\right)+b\sin\left(\frac{\gf}{p-1}\right)
$$
for some real numbers $a,b$. The boundary conditions are the following
\bel{1-2-8}\BA{lll}
(i)\qquad -\myfrac{b}{p-1}+|a|^{p-1}a=0\Longleftrightarrow b=(p-1)|a|^{p-1}a,\\[4mm]
(ii)\qquad-\myfrac{a}{p-1}\sin\left(\myfrac{\gp}{p-1}\right)+\myfrac{b}{p-1}\cos\left(\myfrac{\gp}{p-1}\right)\\[4mm]
+\left(a\cos\left(\myfrac{\gp}{p-1}\right)+b\sin\left(\myfrac{\gp}{p-1}\right)\right)
\left|a\cos\left(\myfrac{\gp}{p-1}\right)+b\sin\left(\myfrac{\gp}{p-1}\right)\right|^{p-1}=0.
\EA\ee
\bth{N2} If $N=2$ the set $\CE$ is always discrete and more precisely,

\noindent1- If $\frac{1}{p-1}\in \BBN^*$, then $0$ is the unique solution to $(\ref{1-2-7})$.

\noindent 2- If $\frac{1}{p-1}\notin \BBN^*$, then $(\ref{1-2-8})$ admits  three solutions $\gw_s$, $-\gw_s$ and zero. 
Furthermore $\gw_s$ keeps a constant sign if $p\geq 2$.
\es
\Proof Because of $(\ref{1-2-8})$-(i) we can assume $a,b>0$. Set $X=(p-1)a^{p-1}$ and 
$$\BA {lll}\Gf(X)=
-\sin\left(\myfrac{\gp}{p-1}\right)+X\cos\left(\myfrac{\gp}{p-1}\right)\\[4mm]
\phantom{\Gf(X)}
+X\left(\cos\left(\myfrac{\gp}{p-1}\right)+X\sin\left(\myfrac{\gp}{p-1}\right)\right)
\left|\cos\left(\myfrac{\gp}{p-1}\right)+X\sin\left(\myfrac{\gp}{p-1}\right)\right|^{p-1}\!\!\!.
\EA$$
All the separable solutions with $a>0$ (and similarly with $a<0$) are obtained with $a^{p-1}=X_0$ and $b=\left(\frac a{p-1}\right)^{\frac 1p}$ where $X_0$ is a positive zero of the function $\Gf$.  \smallskip

\noindent (1) If $\frac{\gp}{p-1}=\frac\gp2+k\gp\;\text{for some } k\in\BBN,$
then $\Gf(X)=(-1)^{k+1}(1-X^{p+1})$. Hence there exist only three solutions corresponding to 
$$(a,b)=(0,0)$$
$$(a,b)=\left(\left(\frac{1}{p-1}\right)^{\frac{1}{p-1}}, \left(\frac{1}{p-1}\right)^{\frac{1}{p-1}}\right)$$
$$ (a,b)=\left(-\left(\frac{1}{p-1}\right)^{\frac{1}{p-1}}, -\left(\frac{1}{p-1}\right)^{\frac{1}{p-1}}\right).
$$

\noindent (2) If $\frac{\gp}{p-1}=k\gp\;\text{for some } k\in\BBN,$
then $\Gf(X)=2(-1)^{k}X$. Hence the only solution is $(a,b)=(0,0)$.\smallskip

\noindent (3) If $\frac{\gp}{p-1}\neq\frac{k\gp}{2}\;\text{for any } k\in\BBN,$ then
$$\BA {lll}\Gf(X)=X-\tan\left(\myfrac{\gp}{p-1}\right)\\[4mm]\phantom{}
+X\tan\left(\myfrac{\gp}{p-1}\right)
\left|\sin\left(\myfrac{\gp}{p-1}\right)\right|^{p-1}\left|\cot\left(\myfrac{\gp}{p-1}\right)+X\right|^{p-1}\left(\cot\left(\myfrac{\gp}{p-1}\right)+X\right).
\EA$$
Hence
$$\BA {lll}\Gf'(X)=1\\[4mm]
+\tan\left(\myfrac{\gp}{p-1}\right)
\left|\sin\left(\myfrac{\gp}{p-1}\right)\right|^{p-1}\left|\cot\left(\myfrac{\gp}{p-1}\right)+X\right|^{p-1}
\left(\cot\left(\myfrac{\gp}{p-1}\right)+(p+1)X\right).
\EA$$
If $\tan\left(\frac{\gp}{p-1}\right)>0$, then $\Gf'(X)>0$ and since $\Gf(0)<0$, $\Gf$ admits a unique root $X_0>0$. Hence 
there exist only three solutions, $\gw_s,-\gw_s,0$. \smallskip

\noindent If $\tan\left(\frac{\gp}{p-1}\right)<0$, then $\Gf(0)>0$. Moreover 
$$\BA{lll}
\Gf''(X)=p\tan\left(\myfrac{\gp}{p-1}\right)\left|\sin\left(\myfrac{\gp}{p-1}\right)\right|^{p-1}\left|\cot\left(\myfrac{\gp}{p-1}\right)+X\right|^{p-3}\\[4mm]\phantom{-------}
{\ti} \left(\cot\left(\myfrac{\gp}{p-1}\right)+X\right)\left((p+1)X+2\cot\left(\myfrac{\gp}{p-1}\right)\right).
\EA$$
Hence $\Gf''$ is negative in the interval $(0,-\frac{2}{p+1}\cot\left(\frac{\gp}{p-1}\right))$, positive in  \\
$(-\frac{2}{p+1}\cot\left(\frac{\gp}{p-1}\right),-\cot\left(\frac{\gp}{p-1}\right))$ and negative in $(-\cot\left(\frac{\gp}{p-1}\right),\infty)$. A standard study shows that $\Gf'$ is positive on $(0, X_*)$ for some 
$X_*>-\cot\left(\frac{\gp}{p-1}\right)$, vanishes at $X_*$ and is negative on $(X_*,\infty)$. Finally, $\Gf$ is increasing on 
$(-\infty,X_*)$ with a positive maximum and negative on $(X_*,\infty)$. As a consequence $\Gf$ admits a unique zero at 
$X_0>-\cot\left(\frac{\gp}{p-1}\right)$ and there exist again only three solutions $\gw_s,-\gw_s$ and $0$. This ends the proof.
\qeda
\mysection{Isolated singularities}
\subsection{Regularity results}
We assume that $\Gw\subset\BBR^N$, $N\geq 2$ is a bounded smooth domain such that $0\in\prt\Gw$. We have the following basic estimate the proof of which is based upon Moser's iterative scheme. 
\bprop {a-priori}Let $g:\BBR\mapsto\BBR$ be a continuous function such that $rg(r)\geq 0$ on $\BBR$. Then any function $u\in C^1(\overline\Gw\setminus\{0\})$ which verifies 
\bel{I-3-1}\BA {lll}
\phantom{g(u;}
-\Gd u=0\qquad&\text{in }\;\Gw\\[0mm]
\myfrac{\prt u}{\prt {\bf n}}+g(u)=0\qquad&\text{in }\;\prt\Gw\setminus\{0\},
\EA
\ee
satisfies for any $a>1$ and some $c_a>0$,
\bel{I-3-2}
\norm u_{L^{\infty}(\Gw\cap B^c_{2r})}\leq \myfrac{c_a}{r^{\frac Na}}\norm u_{L^{a}(\Gw\cap B^c_{r})}\qquad\text{for all }r\in (0,r_0],
\ee
where $r_0>0$ depends on $\Gw$. In particular, if $u$ is nonnegative, then for any $\ge>0$ there exists $c_\ge>0$ such that
\bel{I-3-3}
\norm u_{L^{\infty}(\Gw\cap B^c_{2r})}\leq \myfrac{c_\ge}{r^{N-1+\ge}}\myint{\prt\Gw}{}d\gl,
\ee
where $\gl\in\mathfrak M_+(\prt\Gw)$ is the boundary trace of $u$.
\es

\noindent\Proof Let $\gz\in C^{1,1}(\BBR^N)$ such that $0\leq\gz\leq 1$, $\gz(x)=1$ if $|x|\geq s$, $\gz(x)=0$ if $|x|\leq r$ for some $0<r<s$, and 
$|\nabla\gz(x)|\leq \frac{2}{s-r}\chi_{_{\Gg_r^s}}(x)$ where $\Gg_r^s=\{x\in\BBR^N:r\leq|x|\leq s\}$. For $\ga>0$ we have from $(\ref{I-3-1})$, 
$$\myint{\Gw}{}\langle\nabla u,\nabla(\gz^2|u|^{\ga-1}u )\rangle dx+\myint{\prt\Gw}{}\gz^2g(u)|u|^{\ga-1}u\, dS=0.
$$
Then
$$\BA {lll}\myint{\Gw}{}\langle\nabla u.\nabla(\gz^2|u|^{\ga-1}u )\rangle dx=\myfrac{4\ga}{(\ga+1)^2}\myint{\Gw}{}\left|\nabla |u|^{\frac{\ga+1}{2}}\right|^2\gz^2 dx\\[4mm]
\phantom{\myint{\Gw}{}\nabla u.\nabla(\gz^2|u|^{\ga-1}u )dx}
+\myfrac{4\ga}{\ga+1}\myint{\Gw}{}\gz|u|^{\frac{\ga-1}{2}}u\nabla |u|^{\frac{\ga+1}{2}}.\nabla\gz dx\\[4mm]
\phantom{\myint{\Gw}{}\nabla u.\nabla(\gz^2|u|^{\ga-1}u )dx}
\geq \myfrac{4\ga}{(\ga+1)^2}\myint{\Gw}{}\left|\nabla |u|^{\frac{\ga+1}{2}}\right|^2\gz^2 dx\\[4mm]
\phantom{\myint{\Gw}{}\nabla u.\nabla(\gz^2|u|^{\ga-1}u )dx}
-\myfrac{4\ga}{\ga+1}\left(\myint{\Gw}{}|u|^{\ga+1}|\nabla\gz|^2 dx\right)^{\frac 12}\left(\myint{\Gw}{}\left|\nabla |u|^{\frac{\ga+1}{2}}\right|^2\gz^2dx\right)^{\frac 12}
\EA$$
Put 
$$X=\left(\myint{\Gw}{}\left|\nabla |u|^{\frac{\ga+1}{2}}\right|^2\gz^2dx\right)^{\frac 12},\;
Y=\left(\myint{\Gw}{}|u|^{\ga+1}|\nabla\gz|^2 dx\right)^{\frac 12}
$$
and $A=\myint{\prt\Gw}{}\gz^2g(u)|u|^{\ga-1}u \,dS$, then 
\bel{I-3-4}4\ga X^2-4(\ga+1)XY+(\ga+1)^2A^2\leq 0\Longrightarrow X\leq\myfrac{\ga+1}{\ga}Y.
\ee
The discriminant of this equation in $X$ is necessarily nonnegative, therefore 
\bel{I-3-5}Y\leq\ga A^2\Longleftrightarrow \myint{\prt\Gw}{}\gz^2g(u)|u|^{\ga-1}u \,dS\leq\myfrac{1}{\ga}\myint{\Gw}{}|u|^{\ga+1}|\nabla\gz|^2 dx.
\ee
Since $\gz\nabla |u|^{\frac{\ga+1}{2}}=\nabla \left(\gz |u|^{\frac{\ga+1}{2}}\right)-|u|^{\frac{\ga+1}{2}}\nabla \gz$, we deduce from $(\ref{I-3-4})$ with the help of Young's inequality, 
$$
\myfrac{\ga}{(\ga+1)^2}\myint{\Gw}{}\left|\nabla\left(\gz |u|^{\frac{\ga+1}{2}}\right)\right|^2dx+
\myint{\prt\Gw}{}\gz^2g(u)|u|^{\ga-1}u \,dS\leq \myfrac{13}{\ga}\myint{\Gw}{}|u|^{\ga+1}|\nabla\gz|^2 dx,
$$
which leads to
\bel{I-3-6}\BA {lll}
\myfrac{\ga}{(\ga+1)^2}\norm{\gz |u|^{\frac{\ga+1}{2}}}^2_{W^{1,2}(\Gw)}+\myint{\prt\Gw}{}\gz^2g(u)|u|^{\ga-1}u \,dS\\[4mm]
\phantom{--------}\leq \myfrac{13}{\ga}\myint{\Gw}{}|u|^{\ga+1}|\nabla\gz|^2 dx
+\myfrac{\ga}{(\ga+1)^2}\norm{\gz |u|^{\frac{\ga+1}{2}}}^2_{L^{2}(\Gw)}.
\EA\ee
We first assume $N\geq 3$ and set $\gth=\frac{N}{N-2}$. If $s-r\leq 1$, we obtain, using Gagliardo-Nirenberg's inequality, 
\bel{I-3-7}\BA {lll}
\myfrac{\ga}{(\ga+1)^2}\norm{u}^{\ga+1}_{L^{\gth(\ga+1)}(\Gw\cap B^c_s)}+\myint{\prt\Gw\cap B^c_s}{}g(u)|u|^{\ga-1}u \,dS\leq \myfrac{c^2_{_N}}{\ga(s-r)^2}\norm{u}^{\ga+1}_{L^{\ga+1}(\Gw\cap B^c_r)}.
\EA\ee
We fix $r>0$ and define the sequences for $n\in\BBN^*$
$$\BA{lll}
p_n=\gth p_{n-1}\qquad&\text{with }p_0=\ga+1=a>1\\
r_n=r(2-2^{-n})\qquad&\text{with }r_0=r\\
s_n=r(2-2^{-n-1}),
\EA$$
thus $s_n-r_n=2^{-n-1}r$. We obtain from $(\ref{I-3-7})$
\bel{I-3-8}\BA {lll}
\myfrac{p_n-1}{p_n^2}\norm{u}^{p_n}_{L^{p_{n+1}}(\Gw\cap B^c_{s_n})}+\myint{\prt\Gw\cap B^c_{s_n}}{}g(u)|u|^{p_n-2}u \,dS\leq \myfrac{c_{_N}}{(p_n-1)(s_n-r_n)^2}\norm{u}^{p_n}_{L^{p_n}(\Gw\cap B^c_{r_n})}.
\EA\ee
Therefore
\bel{I-3-9}\BA {lll}
\norm{u}_{L^{p_{n+1}}(\Gw\cap B^c_{s_n})}\leq \left(\myfrac{c_{_N}2^{n+1}(\ga+1)}{\ga r}\right)^{\frac{2}{p_n}}\norm{u}_{L^{p_n}(\Gw\cap B^c_{r_n})}
\EA\ee
Because $s_n\to 2r$ when $n\to\infty$, we obtain by an easy induction
\bel{I-3-10}\BA {lll}
\norm{u}_{L^{\infty}(\Gw\cap B^c_{2r})}\leq \myfrac{c_{N,\ga}}{r^{\frac{N}{a}}}\norm{u}_{L^{a}(\Gw\cap B^c_r)}.
\EA\ee
We notice that we have neglected the boundary integral in $(\ref{I-3-8})$. Indeed, the same induction yields
\bel{I-3-11}\BA {lll}
\norm{u}_{L^{\infty}(\prt\Gw\cap B^c_{2r})}\leq \myfrac{c'_{N,\ga}}{r^{\frac{N}{a}}}\norm{u}_{L^{a}(\Gw\cap B^c_r)}.
\EA\ee

\noindent If $N=2$ we use the interpolation inequality
\bel{I-3-11a}\BA {lll}
\norm{\gz |u|^{\frac{\ga+1}{2}}}^2_{W^{\frac12,2}(\Gw)}\leq c_1\norm{\gz |u|^{\frac{\ga+1}{2}}}_{W^{1,2}(\Gw)}\norm{\gz |u|^{\frac{\ga+1}{2}}}_{L^{2}(\Gw)}.
\EA\ee
Combining it with the imbedding inequality
\bel{I-3-11b}\BA {lll}
\norm{\gz |u|^{\frac{\ga+1}{2}}}^2_{L^{\frac{2N}{N-1}}(\Gw)}\leq c_2\norm{\gz |u|^{\frac{\ga+1}{2}}}^2_{W^{\frac12,2}(\Gw)},
\EA\ee
we obtained that $(\ref{I-3-7})$ is replaced by
\bel{I-3-11c}\BA {lll}
\myfrac{\ga}{(\ga+1)^2}\norm{u}^{\ga+1}_{L^{\tilde\gth(\ga+1)}(\Gw\cap B^c_s)}+\myfrac{1}{2}\myint{\prt\Gw\cap B^c_s}{}g(u)|u|^{\ga-1}u \,dS\leq \myfrac{\tilde c^2_{_N}}{\ga(s-r)^2}\norm{u}^{\ga+1}_{L^{\ga+1}(\Gw\cap B^c_r)}.
\EA\ee
with $\tilde\gth=\frac{N}{N-1}$. {\it Mutatis mutandis}, the end of the proof follows easily.

Next we assume that $u\geq 0$. Then it admits a boundary trace (see e.g. \cite{MaVe-book}) which is a nonnegative Radon $\gl$ measure on $\prt\Gw$ and the Riesz-Herglotz representation formula in terms of Poisson potential of the measure $\gl$ holds,
\bel{I-3-12}\BA {lll}
u(x)={\bf P}^{\Gw}[\gl]:=\myint{\prt\Gw}{}P^\Gw(x,y)d\gl(y)\qquad\text{for all }x\in\Gw,
\EA\ee
where $P^\Gw$ is the Poisson kernel defined in $\Gw\ti\prt\Gw$. Furthermore $u$ belongs to the Lorentz space $L^{\frac{N}{N-1},\infty}(\Gw)$ and $L^{\frac{N+1}{N-1},\infty}_\gr(\Gw)$, where $\gr(x)=\dist (x,\prt\Gw)$ (see e.g. \cite{GmVe}). Furthermore
\bel{I-3-13}\BA {lll}
\norm{u}_{L^{\frac{N}{N-1},\infty}(\Gw)}+\norm{u}_{L^{\frac{N+1}{N-1},\infty}_\gr(\Gw)}\leq c_{\Gw}\norm\gl_{\mathfrak M(\prt\Gw)}. 
\EA\ee
For any $\ge>0$ there exists $c_\ge>0$ such that 
$$\norm{u}_{L^{\frac{N}{N-1+\ge},\infty}(\Gw)}\leq c_\ge\norm{u}_{L^{\frac{N}{N-1},\infty}(\Gw)}
$$
If we apply $(\ref{I-3-11})$ with $a=\frac{N}{N-1+\ge}$ we infer
\bel{I-3-14}\BA {lll}
\norm{u}_{L^{\infty}(\Gw\cap B^c_{2r})}\leq \myfrac{c'_\ge}{r^{N-1+\ge}}\norm\gl_{\mathfrak M(\prt\Gw)}.
\EA\ee
This ends the proof. \qeda
\medskip

\noindent\Remark A natural question is whether  $(\ref{I-3-3})$ is valid with $\ge=0$. Notice that using the standard estimates on the Poisson kernel we have, 
\bel{I-3-15}\BA {lll}
\norm{u}_{L^{\infty}(\Gw_r)}\leq \myfrac{c}{r^{N-1}}\norm\gl_{\mathfrak M(\prt\Gw)}\quad\text{for all }r>0,
\EA\ee
where  $\Gw_r=\{x\in\Gw:\gr(x)\geq r\}$ and $c=c(\Gw)>0$. 
\subsection{Linear estimates}
We assume that $\Gw$ is a bounded smooth domain of $\BBR^N$, $N\geq 2$.
\bprop{lin} Let $a\geq 0$ be a constant and $\gl$ and $\gm$ be two bounded Radon measures on $\Gw$ and $\prt\Gw$ respectively. Then there exists a unique weak solution 
$u\in L^1(\Gw)$ of
\bel{I-L-1}\BA{lll}
-\Gd u+au =\gl\qquad&\text{in }\Gw\\[0mm]\phantom{---,}
\myfrac{\prt u}{\prt {\bf n}}=\gm\qquad&\text{in }\;\Gw.
\EA
\ee
Furthermore there exists $c=c(\Gw)>0$ such that 
\bel{I-L-2}\BA{lll}
\norm u_{L^{\frac{N}{N-2},\infty}(\Gw)}+\norm {\nabla u}_{L^{\frac{N}{N-1},\infty}(\Gw)}\leq c\left(\norm \gl_{\mathfrak M(\Gw)}+\norm \gm_{\mathfrak M(\prt\Gw)}\right) +b_a\norm u_{L^1(\Gw)},
\EA
\ee
if $N>2$ with $b_a\geq 0$, $b_a>0$ if $a=0$,  and 
\bel{I-L-3}\BA{lll}
\norm u_{L^{r}(\Gw)}+\norm {\nabla u}_{L^{2,\infty}(\Gw)}\leq c(r)\left(\norm \gl_{\mathfrak M(\Gw)}+\norm \gm_{\mathfrak M(\prt\Gw)}\right)+b_a\norm u_{L^1(\Gw)},
\EA
\ee
for any $r<\infty$, if $N=2$. 
\es
\Proof We first consider the case $\Gw=B_R^+:=\{x=(x',x_{_N})\in B_R:x_{_N}>0\}$. We set
$$\tilde u(x',x_{_N})=\left\{\BA{lll} u(x',x_{_N})\qquad&\text{if }x_{_N}>0\\[1mm]
u(x',-x_{_N})\qquad&\text{if }x_{_N}<0.
\EA\right.
$$
Then $\tilde u$ satisfies
\bel{I-L-4}\BA{lll}
-\Gd \tilde u+a\tilde u =\tilde\gl+2\CH_{\prt\BBR^N_+}\gm\qquad&\text{in }B_R\\[0mm]\phantom{---,}
\myfrac{\prt \tilde u}{\prt {\bf n}}=\tilde \gm\lfloor_{\prt B_R}\qquad&\text{in }\;\prt B_R,
\EA
\ee
where $\CH_{_{\prt\BBR^N_+}}$ is the (N-1)-dimensional Hausdorff measure and $\tilde\gl$ and $\tilde\gm$ are defined accordingly to $\tilde u$ by an even reflexion through $\prt\BBR^N_+$. Then $\tilde u$ satisfies locally $(\ref{I-L-2})$ in the sense that for any $0<R'<R$ there holds
\bel{I-L-5}\BA{lll}
\norm {\tilde u}_{L^{\frac{N}{N-2},\infty}(B_{R'})}+\norm {\nabla \tilde u}_{L^{\frac{N}{N-1},\infty}((B_{R'})}
\leq c\left(\norm {\tilde \gl}_{\mathfrak M(B_R)}+\norm \gm_{\mathfrak M(\prt B_{R'}^+)}\right),
\EA
\ee
when $N>2$, with straightforward modification if $N=2$. This implies 
\bel{I-L-6}\BA{lll}
\norm  u_{L^{\frac{N}{N-2},\infty}(B^+_{R'})}+\norm {\nabla  u}_{L^{\frac{N}{N-1},\infty}((B^+_{R'})}\leq c\left(\norm {\tilde \gl}_{\mathfrak M(B_R)}+\norm {\gm}_{\mathfrak M(\prt B^+_R)}\right);
\EA
\ee
For a general domain $\Gw$, consider a point $a\in\prt\Gw$. There exists $r_a>0$ such that we can perform an even reflexion though $\prt\Gw\cap B_{r_a}(a)$  following the normal vector to $\prt\Gw$ as in \cite [Lemma 2.4]{BorVe}, with the modification that we use an even reflection and not the odd one which is therein adapted to zero boundary data. If we denote by $\tilde u$ the reflected function defined in $B_{r_a}(a)$, it satisfies
\bel{I-L-7}\BA{lll}\displaystyle
-\sum_j\myfrac{\prt}{\prt x_j}\tilde A_j(x,\nabla\tilde u)+a\tilde u =\tilde\gl+2\CH_{\prt\BBR^N_+}\gm\qquad&\text{in }B_{r_a}(a),
\EA
\ee
where the $A_j$ are $C^1$ functions satisfying the standard ellipticity and boundedness conditions. The local regularity theory yields 
\bel{I-L-8}\BA{lll}
\norm {\tilde u}_{L^{\frac{N}{N-2},\infty}(B_{r_a'}(a))}\!\!+\norm {\nabla \tilde u}_{L^{\frac{N}{N-1},\infty}((B_{r_a'}(a))}
\!\!\leq c\left(\norm {\tilde \gl}_{\mathfrak M (B_{r_a}(a))}\!\!+\norm \gm_{\mathfrak M(\prt B^+_{{r_a}}(a))}\right)\!.
\EA
\ee
for any $0<r'_a<r_a$, where $c$ depends on $\Gw$ and $r_a-r'_a$. We obtain $(\ref{I-L-2})$ by a compactness argument. The proof of $(\ref{I-L-3})$ is similar.
Uniqueness is straightforward. \qeda\medskip

\noindent \Remark These results are not new. However they show that the estimates are local which will be useful later on in the sense that for any compact set 
$K\subset\overline\Gw$ and any $\ge>0$ there holds
 \bel{I-L-8'}\BA{lll}
\norm {u}_{L^{\frac{N}{N-2},\infty}(K)}\!\!+\norm {\nabla \tilde u}_{L^{\frac{N}{N-1},\infty}(K)}
\!\!\leq c\left(\norm {\tilde \gl}_{\mathfrak M (\Gw\cap K_\ge)}\!\!+\norm \gm_{\mathfrak M(\prt \Gw\cap K_\ge)}\right)\!,
\EA
\ee
where $K_\ge=\{x\in\BBR^N:\dist(x,K)\leq\ge\}$ and where $c$ is a positive constant depending on $\Gw,K,\ge$.\medskip

\noindent \Remark A more general global statement of existence and regularity with a more involved proof can be found in \cite[Theorems 1, 2]{MeRa}. The same estimates holds up to replacing $b_a\norm u_{L^1(\Gw)}$ by $b_a\norm u_{L^1(\prt\Gw)}$ in $(\ref{I-L-2})$-$(\ref{I-L-3})$  if $(\ref{I-L-1})$
is replaced by
\bel{I-L-1'}\BA{lll}
\phantom{,,,}-\Gd u =\gl\qquad&\text{in }\Gw\\[0mm]
\myfrac{\prt u}{\prt {\bf n}}+au=\gm\qquad&\text{in }\;\Gw.
\EA
\ee

\blemma {lem-br} Let $\gl\in L^{1}(\Gw)$, $\gm \in L^{1}(\prt\Gw)$ and $u\in L^1(\Gw)$ be the weak solution of $(\ref{I-L-1})$. Then we have 
 for all $\gz\in \CC(\Gw)$, $\gz\geq 0$,
\bel{I-L-9}\BA{lll}
\myint{\Gw}{}|u|(-\Gd \gz+a\gz)dx\leq \myint{\Gw}{}\gl\gz {\rm sign}_0 (u) dx+\myint{\prt\Gw}{}\gm\gz{\rm sign}_0(u) dS
\EA
\ee 
where ${\rm sign}_0 =\chi_{_{(0,\infty)}}-\chi_{_{(-\infty,0)}}$,  and
\bel{I-L-10}\BA{lll}
\myint{\Gw}{}u_+(-\Gd \gz+a\gz)dx\leq \myint{\Gw}{}\gl\gz {\rm sign}_0^+ (u) dx+\myint{\prt\Gw}{}\gm\gz{\rm sign}^+_0(u) dS
\EA
\ee 
where ${\rm sign}_0^+=\chi_{_{(0,\infty)}}$.
\es
\noindent\Proof We first assume that $u$ is a smooth function. Let $\{\gg_k\}\subset C^\infty_0(\BBR)$ be a sequence of nonnegative functions such that $0\leq\gg_k\leq 1$, $\gg_k=0$ on $(-\infty,0]$, $\gg_k'\geq 0$, $\gg_k=1$  on $[k^{-1},\infty)$, and let $\gz\in \CC(\Gw)$, $\gz\geq 0$. Then 
$$\BA {lll}\myint{\Gw}{}\left(\gg_k'(u)|\nabla u|^2\gz +\gg_k(u)\langle\nabla u,\nabla\gz\rangle \right)dx+a\myint{\Gw}{}u\gg_k(u) dx\\[4mm]
\phantom{--------------------}=\myint{\Gw}{}\gg_k(u)\gz\gl dx+\myint{\prt\Gw}{}\gg_k(u)\gz\gm dS.
\EA$$
 Set $j_k(r)=\int_0^r\gg_k(s) ds$, then
 $$\myint{\Gw}{}\langle\nabla j_k(u),\nabla\gz\rangle dx+a\myint{\Gw}{}u\gg_k(u) dx\leq \myint{\Gw}{}\gg_k(u)\gz\gl dx+\myint{\prt\Gw}{}\gg_k(u)\gz\gm dS.
$$
Since $\gz\in\CC(\Gw)$, 
\bel{I-L-11}\myint{\Gw}{}\left(-j_k(u)\Gd \gz +au\gg_k(u)\right)dx\leq \myint{\Gw}{}\gg_k(u)\gz\gl dx+\myint{\prt\Gw}{}\gg_k(u)\gz\gm dS.
\ee
Letting $k\to\infty$, we infer
$$\myint{\Gw}{}\left(-\Gd \gz +a\gz\right)u_+dx\leq \myint{\Gw}{}sign_+ (u)\gz\gl dx+\myint{\prt\Gw}{}sign_+ (u)\gz\gm dS.
$$
In the same way, we prove
$$\myint{\Gw}{}\left(-\Gd \gz +a\gz\right)|u|dx\leq \myint{\Gw}{}sign (u)\gz\gl dx+\myint{\prt\Gw}{}sign (u)\gz\gm dS.
$$
In the general case, let $\{\gl_\ell\}$, $\{\gm_\ell\}$ be two sequences converging in $L^1(\Gw)$ and $L^{1}(\prt\Gw)$ to $\gl$ and $\gm$ respectively. Then the sequence of solutions $\{u_\ell\}$ of 
\bel{I-L-1n}\BA{lll}
-\Gd u_\ell+au_\ell =\gl_\ell \qquad&\text{in }\Gw\\[0mm]\phantom{_\ell---,}
\myfrac{\prt u_\ell}{\prt {\bf n}}=\gm_\ell \qquad&\text{in }\;\Gw,
\EA
\ee
converges to the solution $u$ of $(\ref{I-L-1n})$ in $L^{\frac{N}{N-2},\infty}(\Gw)$ (any $L^r(\Gw)$ with $1<r<\infty$ if $N=2$) and $\{\nabla u_\ell\}$ converges to 
$\nabla u$ in $\left(L^{\frac{N}{N-1},\infty}(\Gw)\right)^N$. This implies that $(\ref{I-L-11})$ holds, hence $(\ref{I-L-9})$ and $(\ref{I-L-10})$ follow. \qeda
\medskip

The following general regularity result proved in \cite[Theorem 6]{QR} will be used later on.

\bprop{QR} Let $a\geq 0$ be a constant and $u$ be the weak solution of $(\ref{I-L-1})$ with $\gl=0$ and $\gm\in L^{m}(\prt\Gw)$, $m\geq 1$. Let 
$q\in [m,\infty]$. The following regularity results hold:\smallskip

\noindent 1- If $\frac{1}{m}-\frac{1}{q}<\frac{1}{N-1}$, then $u\in L^q(\prt\Gw)$.\smallskip

\noindent 2- If $\frac{1}{m}-\frac{N}{(N-1)q}<\frac{1}{N-1}$, then $u\in L^q(\Gw)$.\smallskip

\noindent 3- If $\frac{1}{m}-\frac{N}{(N-1)q}<0$, then $u\in W^{1,q}(\Gw)$.
\es 

\noindent\Remark In each case of the above proposition there holds
\bel{I-L-12}
\norm u_{\BBX}\leq c\norm \gm_{L^{m}(\prt\Gw)}+b_a\norm u_{L^1(\Gw)},
\ee
where $\BBX$ is either $L^q(\prt\Gw)$ either $L^q(\Gw)$ or $W^{1,q}(\Gw)$ and $b_a\geq 0$ is as in \rprop{lin}. 
From this result we obtain higher regularity according to the regularity of the boundary data.

\bprop{Reg} Let $a\geq 0$ be a constant and $u$ be the weak solution of $(\ref{I-L-1})$ with $\gl=0$. If $\gm\in W^{1,m}(\prt\Gw)$, $m\geq 1$. Let 
$q\in [m,\infty]$ be such that $\frac{1}{m}-\frac{N}{(N-1)q}<0$. Then $u\in W^{2,q}(\Gw)$. Moreover
\bel{I-L-13}
\norm u_{W^{2,q}(\Gw)}\leq c\norm \gm_{W^{1,m}(\prt\Gw)}+b_a\norm u_{L^1(\Gw)}.
\ee
\es 
\Proof For the sake of simplicity we assume that $\Gw=B_1$, the unit ball in $\BBR^N$. In spherical coordinates $u$
 satisfies
 \bel{I-L-14}\BA {lll}
-u_{rr}-\myfrac{N-1}{r}u_r-\myfrac{1}{r^2}\Gd'u=0\qquad&\text{in }(0,1)\ti S^{N-1}\\[2mm]
\phantom{-,--u_r----}
u_r(1,.)=\gm(.)\qquad&\text{in } S^{N-1},
\EA\ee
where $\Gd'$ is the Laplace-Beltrami operator on $S^{N-1}$. Let $A$ be a skew-symmetric matrix in $\BBR^N$, $X_t:=\exp(tA)$ the group of isometries that it generates and $L_A$ the Lie derivative defined by
$$L_Aw(\gs)=\myfrac{d}{dt}w(X_t\gs)\lfloor_{t=0}
$$
Since $L_A$ commutes with $\Gd'$, the function $(r,\gs)\mapsto v(r,\gs)=L_Au(r,\gs)$ satisfies
$$\BA {lll}
-v_{rr}-\myfrac{N-1}{r}v_r-\myfrac{1}{r^2}\Gd'v=0\qquad&\text{in }(0,1)\ti S^{N-1}\\[2mm]
\phantom{-,--v_r----}
v_r(1,.)=L_A\gm(.)\qquad&\text{in } S^{N-1}.
\EA$$
We deduce
$$
\norm v_{W^{1,q}(\Gw)}\leq c\norm {L_A\gm}_{L^{m}(\prt\Gw)}+b_a\norm u_{L^1(\Gw)}\leq c\norm \gm_{W^{1,m}(\prt\Gw)}+b_a\norm u_{L^1(\Gw)}.
$$
This implies firstly that 
$$\norm {\nabla'v}_{W^{1,q}(\Gw)}\leq c'\norm \gm_{W^{1,m}(\prt\Gw)}+Nb_a\norm u_{L^1(\Gw)},
$$
which is an estimate for all the tangential derivatives of $v$ and we obtained the final estimate with the normal derivative using the equation. \qeda
\medskip

Interating this method and using interpolation techniques, we obtain
\bprop{highReg} Let $a\geq 0$ be a constant and $u$ be the weak solution of $(\ref{I-L-1})$. If $\gm\in W^{k+s,m}(\prt\Gw)$, $m\geq 1$, $k\in\BBN^*$, $s\in (0,1)$. Let 
$q\in [m,\infty]$ be such that $\frac{1}{m}-\frac{N}{(N-1)q}<0$. If $\gl\in W^{k-1+s,q}(\Gw)$ then $u\in W^{k+1+s,q}(\Gw)$ and
\bel{I-L-13a}
\norm u_{W^{k+1+s,q}(\Gw)}\leq c\left(\norm \gm_{W^{k+s,m}(\prt\Gw)}+\norm \gl_{W^{k-1+s,q}(\Gw)}\right)+b_a\norm u_{L^1(\Gw)}.
\ee
\es 

The next local version of the previous results will be used later on.
\bprop{loc} Let $a\geq 0$ be a constant, $N\subsetneq\prt\Gw$ be compact and $u$ be a nonnegative weak solution of $(\ref{I-L-1})$ with $\gl=0$ and  $\gm\in W_{loc}^{k+s,m}(\prt\Gw\setminus N)$, $m\geq 1$, $k\in\BBN^*$, $s\in (0,1)$. Let 
$q\in [m,\infty]$ be such that $\frac{1}{m}-\frac{N}{(N-1)q}<0$. If for $\ge>0$ small enough we set $N_\ge=\{x\in\overline\Gw:\dist (x,N)\leq\ge\}$, then $u\in W^{k+1+s,q}(\Gw\setminus N_\ge)$ and
\bel{I-L-13b}
\norm u_{W^{k+1+s,q}(\Gw\setminus N_{2\ge})}\leq c\norm \gm_{W^{k+s,m}(\prt\Gw\setminus N_\ge)}+b_a\norm u_{L^1(\Gw\setminus N_\ge)},
\ee
with $c=c(\ge)>0$. 
\es 
\Proof Let $\gz\in C^\infty_0(\BBR^N)$, $\gz\geq 0$, vanishing in a neighborhood of $N$ and $v=\gz u$, then 
\bel{I-L-14a}\BA {lll}
-\Gd v=\gz\gl+u\Gd \gz+2\langle\nabla u,\nabla\gz\rangle&\quad\text{in }\Gw\\
[2mm]\phantom{-}
\myfrac{\prt v}{\prt{\bf n}}=\gz\gm-u\myfrac{\prt \gz}{\prt{\bf n}}&\quad\text{in }\prt\Gw
\EA\ee
Since $u$ is positive harmonic, it belongs to $L^{\frac{m}{2-m}}(\Gw)$ and $|\nabla u|\in L^{m}(\Gw)$ for any $m\in (1,\frac{N}{N-1})$. Then 
$\gz\gl+u\Gd \gz+2\nabla u.\nabla\gz\in L^m(\Gw)$. Combining the trace theorem  and 
Sobolev imbedding theorem, $u\frac{\prt \gz}{\prt{\bf n}}\in W^{\frac{m-1}{m},m}(\prt\Gw)\subset L^{\frac{m(N-1)}{N-m}}(\prt\Gw)$. 
Since the solution $w$ of 
$$\BA {lll}
-\Gd w=\gz\gl+u\Gd \gz+2\langle\nabla u,\nabla\gz\rangle&\quad\text{in }\Gw\\
[2mm]\phantom{-}
\myfrac{\prt w}{\prt{\bf n}}=0&\quad\text{in }\prt\Gw
\EA$$
with zero average belongs to $W^{2,m}(\Gw)\subset W^{1,\frac{mN}{N-m}}$, it follows from \rprop{QR} that $v\in W^{1,q}(\Gw)$ for any $q<\frac{mN}{N-m}$. We iterate this process by setting 
$$m_0=m<\frac{N}{N-1}\;\text{ and } \frac{1}{m_n}=\frac{1}{m_{n-1}}-\frac 1N=\frac{1}{m_{0}}-\frac nN\;\text{ for }n\in\BBN^*.$$
If $n^*$ is the largest integer smallest than $\frac{N}{m_0}$, then $v\in W^{1,m_{n^*\!+1}-\gt}(\Gw)$ for any $\gt>0$, hence 
$v\in W^{1,m_{n^*\!+1}-\gt}(\Gw)\subset W^{s,\infty}(\Gw)$ for some $s\in (0,1)$. By \rprop{highReg}, $v\in W^{1+s,\infty}(\Gw)$. Iterating this method we obtain the claim. $\phantom{----}$\qeda\medskip

\noindent\Remark The sign assumption on $u$ may look unusual, but it must be noticed that the problem is by essence non-local. The only local aspect 
is the one dealing with the local properties of nonnegative harmonic functions and the solutions of elliptic equations with measure data. If we want to get rid of it, we need $\nabla u\in L^{\frac{N}{N-1},\infty}(K)$ for all compact set  $K\subset \overline\Gw\setminus N$ as a starting point of the proof of \rprop{loc}. \medskip

An important application deals with nonlinear boundary value such as
\bel{I-L-15}\BA {lll}
\phantom{-,-}-\Gd u=0&\qquad\text{in }\Gw\\
[1mm]
\myfrac{\prt u}{\prt{\bf n}}+g(u)=0&\qquad\text{in }\prt\Gw,
\EA\ee
where $g:\BBR_+\mapsto\BBR_+$ is a $C^{k+1}$ function. Putting $\gm=g(u)$ and iterating \rprop{loc} we obtain
\bel{I-L-16}\BA {lll}
\norm u_{W^{k+1+s,\infty}(\Gw\setminus N_{2\ge})}\leq c \norm u_{L^\infty(\prt\Gw\setminus N_{\ge})}+b_a\norm u_{L^1(\Gw\setminus N_{\ge})}.
\EA\ee

\blemma{est-sup} Let $1<p\leq \frac{N}{N-1}$ and $u$ be a nonnegative solution of  $(\ref{I-1-9})$ such that $|x|^{\frac{1}{p-1}}u(x)$ is bounded, then $|x|^{\frac{1}{p-1}+\ell}|D^\ell u(x)|$ is also bounded for $\ell=1,2,3$.
\es
\Proof For $k\in (0,1]$ we set $u_k(x)=T_k[u](x)=k^{\frac{1}{p-1}}u(kx)$ where $T_k$ is already defined in $(\ref{I-1-4})$. Then $u_k$ satisfies 
$(\ref{I-1-9})$ in $\Gw^k:=k^{-1}\Gw$. Since $u(x)\leq c|x|^{-\frac{1}{p-1}}$ in $\Gw$, $u_k$ satisfies the same estimate with the same 
constant in $\Gw^k$. Let $r>0 $ such that $\frac{k}{4}\leq r\leq 8k$. By \rprop {loc} we have 
$$\BA {lll}
\norm{u_k}_{W^{3,\infty}(\Gw^k\cap\Gg^{4r/k}_{2r/k})}
\leq c\norm{u_k}_{L^\infty(\prt\Gw^k\cap\Gg_{r/k}^{6r/k})}+b_a\norm{u_k}_{L^1(\Gw^k\cap\Gg_{r/k}^{6r/k})},
\EA$$
where $\Gg_a^b=\{x\in\BBR^N:a\leq|x|\leq b\}$. Since the curvature of $\Gw^k$ is bounded independently of $k$, the constant $c$ is independent of $k$ too. 
Furthermore 
$$\BA {lll}\myint{\Gw^k\cap\Gg_{r/k}^{6r/k}}{}u_k(x) dx= k^{1-N+\frac{1}{p-1}}\myint{\Gw\cap\Gg_{r}^{6r}}{}u(y) dy\\[4mm]
\phantom{\myint{\Gw^k\cap\Gg_{r/k}^{6r/k}}{}u_k(x) dx}\leq c
k^{1-N+\frac{1}{p-1}}\myint{\Gw\cap\Gg_{r}^{6r}}{}|y|^{-\frac{1}{p-1}} dy\\[4mm]
\phantom{\myint{\Gw^k\cap\Gg_{r/k}^{6r/k}}{}u_k(x) dx}=cc_{_N}k^{1-N+\frac{1}{p-1}}\myint{r}{6r}s^{N-2-\frac{1}{p-1}}ds.
\EA$$
This last term is bounded as we have chosen  $\frac{k}{4}\leq r\leq 8k$. Since $D^\ell u_k(x)=k^{\frac{1}{p-1}+\ell}D^ku(kx)$, we take $k=r$ and deduce
$$|D^\ell u(x)|\leq  c'|x|^{-\frac{1}{p-1}-\ell}+c'',
$$
which ends to proof. \qeda

\subsection{Proof of Theorem B}
We denote by $(x,z)\mapsto \CN_\Gw(x,z)$ be the kernel function defined in $\Gw\ti\prt\Gw$ with Neumann boundary data $\gd_z$, that is the solution of 
$v=v_z$ of
\bel{I-B-0}\BA {lll}
-\Gd v+v=0\qquad&\text{in }\Gw\\[0mm]\phantom{---}
\!\myfrac{\prt v}{\prt{\bf n}}=\gd_z&\text{in }\prt\Gw
\EA\ee
It is known that 
\bel{I-B-1}
\CN_\Gw(x,z)\approx\left\{\BA {lll}|x-z|^{2-N}\quad&\text{if }N\geq 3\\[1mm]
-\ln |x-z|&\text{if }N= 2.
\EA
\right.
\ee
Furthermore, if $\gm\in\mathfrak M({\prt\Gw})$ the solution of 
\bel{I-B-0'}\BA {lll}
-\Gd v+v=0\qquad&\text{in }\Gw\\[0mm]\phantom{---}
\!\myfrac{\prt v}{\prt{\bf n}}=\gm&\text{in }\prt\Gw
\EA\ee
is expressed by 
\bel{I-B-0''}\BA {lll}
v(x)=\myint{\prt\Gw}{}\CN_\Gw(x,y)d\gm(y).
\EA\ee

Let $j:\BBR\mapsto\BBR_+$ be a $C^2$ nondecreasing convex function, vanishing on $(-\infty,0]$, such that $0<j'(r)\leq 1$ on $(0,\infty)$. For
 $\ge>0$  set $w_\ge=j(u-\ge \CN_\Gw(.,0))$, then 
 $$-\Gd w_\ge=-\ge \CN_\Gw(.,0)j'(u-\ge \CN_\Gw(.,0))-j''(u-\ge \CN_\Gw(.,0))|\nabla (u-\ge \CN_\Gw(.,0))|^2\leq 0.
 $$
Since $w_\ge$ vanishes in a neighborhood of $0$, 
 $$\myint{\Gw}{}|\nabla w_\ge|^2+\myint{\prt\Gw}{}w_\ge g(u) dS\leq 0. 
 $$
As $g(u)$ has the sign of $u$, it is nonnegative on the support of $w_\ge$. Hence $\nabla w_\ge=0$. This implies that $j(u-\ge \CN_\Gw(.,0))$ is equal to some consatnt $c_\ge$ which is decreasing with $\ge$. Letting $\ge\to 0$ we infer that $u_+$ is constant. Similarly $u_-$ is constant and such is $u$. Notice that 
 for this constant $u$, $g(u)=0$.  $\phantom{--sfzfr--}$ \qeda
\subsection{Proof of Theorem C}
\subsubsection{Straightening the boundary}
If $p>\frac{N-1}{N-2}$, then $u(x)=O(|x|^{-\frac{1}{p-1}})=o(|x|^{2-N})$ and $u=0$ by Theorem B. Therefore we can assume $1<p\leq \frac{N-1}{N-2}$ in the sequel. The basic technique is to straighten the boundary and transform the study near the singular point into a problem in a infinite cylinder. We abridge the proof since the details of the method (initialy introduced in \cite {GmVe} ) can be found in \cite{CV-lin}. We assume that the 
orthonormal basis ${\bf e}_1,...,{\bf e}_{_N}$ is $\BBR^N$ is such that  at $0$, ${\bf n}=-{\bf e}_{_N}$  and that $\prt\Gw$ is locally the graph of a $C^2$ function $\gth$ defined in
$B_{R'}=B_R\cap\{x:x_{_N}=0\}$ and satisfying $\gth(0)=0$, $D\gth(0)=0$. Putting 
$$y_j=x_j=\Gth_j(x)\,\text{ if }j=1,...,N-1\;\text{ and }y_{_N}=x_{_N}-\gth(x')=\Gth_{_N}(x),
$$
then $\Gth=(\Gth_1,...,\Gth_{_N})$ is a local diffeomorphism near $0$. We set $u(x)=\tilde u(y)=\tilde u(r,\gs)=r^{-\frac{1}{p-1}}v(t,\gs)$ with $t=\ln r$. Performing a lengthy computation we derive that $v$ satisfies
\bel{II-C-1}\BA {lll}
(1+\ge_1)v_{tt}+\left(N-\myfrac{2p}{p-1}+\ge_2\right)v_t+\left(\ell_{N,p}+\ge_3\right)v+\Gd'v+\langle\nabla'v,{\ge_4}\rangle\\[2mm]
\phantom{----}+\langle\nabla'v_t,{\ge_5}\rangle
+\langle\nabla'\langle\nabla'v,{\bf e}_{_N}\rangle,{\ge_6}\rangle=0
\EA
\ee
in $(-\infty,T_0]\ti S^{N-1}_+$ and 
\bel{II-C-1'}\BA {lll}
\left(v_t-\myfrac{1}{p-1}v\right)\left(\ge_8-\langle{\bf c},{\bf e}_{_N}\rangle\ge_7\right)-\langle\nabla' v,{\bf e}_{_N}\rangle\ge_7+\langle\nabla'v,{\gn}\rangle+v^p=0
\EA\ee
in $(-\infty,T_0]\ti \prt S^{N-1}_+$, where ${\bf c}=\frac y{|y|}$ and the $\ge_j$ satisfy
\bel{II-C-2}\BA {lll}
|\ge_j(t,.)|+|\ge_{j\,t}(t,.)|+|\nabla'\ge_j(t,.)|\leq ce^t,
\EA\ee
as a consequence of the fact that $|\gth(x')|=0(|x'|^2)$ near $0$. Furthermore the quantities $v_t,v_{tt}, v_{ttt}$ and $\nabla^\ga\prt_{t^\gb}v$ are uniformly bounded on 
$(-\infty,T_1]\ti S^{N-1}_+$ if
$|\ga|+\gb\leq 3$ and $T_1<T_0$ by \rlemma{est-sup}. 
\blemma{energy} There holds
\bel{II-C-3}\BA {lll}
\myint{-\infty}{T_1}\myint{S^{N-1}_+}{}\left(v_t^2+v_{tt}^2+|\nabla'v_t|^2\right) dSdt<\infty.
\EA\ee
\es
\Proof We multiply the first equation $(\ref{II-C-1})$ by $v_t$, integrate on $S^{N-1}_+$ and obtain
\bel{II-C-4}\BA {lll}
\myfrac{1}{2}\myfrac{d}{dt}\left[\myint{S^{N-1}_+}{}\left(v_t^2+\ell_{N,p}v^2-|\nabla'v|^2\right)dS-\myfrac{2}{p+1}
\myint{\prt S^{N-1}_+}{}|v|^{p+1}dS'\right]\\[4mm]\phantom{----------}
+\left(N-\myfrac{2p}{p-1}+\ge_2\right)\myint{S^{N-1}_+}{}v_t^2dS+\eta_1(t)+\eta_2(t)=0,
\EA\ee
where 
$$\eta_1(t)=\myint{S^{N-1}_+}{}\left(\ge_1v_{tt}+\ge_3 v+\langle\nabla 'v,\ge_4\rangle+\langle\nabla' v_t,\ge_5\rangle+\langle\nabla'\langle\nabla' v,{\bf e}_{_N}\rangle,\ge_6\rangle\right) v_tdS,
$$
$$\eta_2(t)=\myint{\prt S^{N-1}_+}{}\!\left(\left(v_t-\myfrac{1}{p-1}v\right)\left(\langle{\bf c},{\bf e}_{_N}\rangle\ge_7-\ge_8\right)+\langle\nabla' v,{\bf e}_{_N}\rangle\ge_7\right) v_tdS'.
$$
By $(\ref{II-C-2})$, $|\eta_j(t)|\leq ce^t$.  The fact that $v_t$ and $\nabla'v$ are uniformly bounded and $N-\frac{2p}{p-1}\neq 0$ as $p\neq \frac{N}{N-2}$, we infer
\bel{II-C-5}\BA {lll}
\myint{-\infty}{T_1}\myint{S^{N-1}_+}{}v_t^2dSdt<\infty.
\EA\ee
Since $vv_{tt}=(vv_t)_t-v_t^2$ and $\langle\nabla'v,\nabla'v_{tt}\rangle=(\langle\nabla'v,\nabla'v_t\rangle)_t-|\nabla'v_t|^2$, we obtain by multiplying $(\ref{II-C-1})$ by $v_{tt}$ and integrating on $S^{N-1}_+$,

\bel{II-C-4'}\BA {lll}
\myfrac{d}{dt}\left[\myint{S^{N-1}_+}{}\left(\left(\myfrac N2-\myfrac{p}{p-1}\right)v_t^2+\ell_{N,p}vv_t-\nabla'v.\nabla'v_t\right)dS
-\myint{\prt S^{N-1}_+}{}|v|^{p-1}vv_tdS'\right]\\[4mm]\phantom{----}
+\myint{S^{N-1}_+}{}\!\!\left(1+\ge_1)(v_{tt}^2-\ell_{N,p}v_t^2+|\nabla'v_t|^2\right)dS+p\myint{\prt S^{N-1}_+}{}|v|^{p-1}v_t^2dS'\\[4mm]
\phantom{----}+\gg_1(t)+\gg_2(t)
=0,
\EA\ee
where 
$$\gg_1(t)=\myint{S^{N-1}_+}{}\left(\ge_2v_t+\ge_3v+\langle\nabla'v,\ge_4\rangle+\langle\nabla'v_t,\ge_5\rangle+\langle\nabla'\langle\nabla'v,{\bf e}_{_N}\rangle,\ge_6\rangle\right)v_{tt}dS,
$$
$$\gg_2(t)=\myint{\prt S^{N-1}_+}{}\left(v_t-\myfrac{1}{p-1}v\right)\left(\langle{\bf c},\nabla'v\rangle\ge_7-\ge_8\right)v_{tt}dS'.
$$
Again $|\gg_j(t)|\leq ce^t$ by $(\ref{II-C-2})$. Since $v_t$, $vv_t$, $\nabla'v.\nabla'v_t$ and $|v|^{p-1}v_t^2$ are uniformly bounded on 
$(-\infty,T_1]\ti S^{N-1}_+$ we infer 
$$\myint{-\infty}{T_1}\myint{S^{N-1}_+}{}\left(v_{tt}^2+|\nabla'v_t|^2\right) dSdt<\infty,$$
which ends the proof.\qeda\medskip

\subsubsection {Strong singularities}
 Because the functions $v_t,v_{tt}$ and $\nabla'v_t$ are uniformly continuous on $(-\infty,T_1]\ti S^{N-1}_+$ we deduce easily from $(\ref {II-C-3})$ that 
\bel{II-C-5'}\lim_{t\to\infty}(v_t^2+v_{tt}^2+|\nabla'v_t|^2)(t,.)=0\quad\text{uniformly on }S^{N-1}_+.
\ee
The negative trajectory of $t\mapsto v(t,.)$ in $C^{2}(S^{N-1}_+)$ is $\CT_-[v]:=\bigcup_{t\leq T_1}\{v(t,.)\}$. By \rlemma{est-sup}, $\CT_-[v]$ is bounded in $C^{3}(S^{N-1}_+)$, hence it is relatively compact in $C^{2}(S^{N-1}_+)$ by the Arzela-Ascoli theorem. Therefore, the alpha-limit set of $\CT_-[v]$ defined by 
\bel{II-C-6}
\CA[\CT_-[v]]:=\bigcap_{t\leq T_1}\text{clo}_{C^2(S^{N-1}_+)}\left(\bigcup_{\gt\leq t}\{v(\gt,.)\}\right)
\ee
is a non-empty compact connected set in $C^2(S^{N-1}_+)$. Using $(\ref{II-C-5})$ and letting $t\to-\infty$ in $(\ref{II-C-1})$, we conclude that if 
$\gw\in \CA[\CT_-[v]]$, then 
\bel{II-C-7}\BA {lll}
\Gd'\gw+\ell_{N,p}\gw=0\qquad&\text{in }S^{N-1}_+\\[2mm]
\phantom{}
\!\!\langle\nabla'\gw,{ \gn}\rangle+\gw^p=0
\qquad&\text{in }\prt S^{N-1}_+,
\EA\ee
hence $\CA[\CT_-[v]]\subset \CE_+$. \smallskip

\noindent If $p=\frac{N-1}{N-2}$, then $\frac{1}{p-1}=N-2$ and $\CE=\{0\}$. Hence $v(t,.)\to 0$ when $t\to-\infty$, equivalently
\bel{II-C-8}\BA {lll}\displaystyle
\lim_{x\to 0}|x|^{N-2}u(x)=0.
\EA\ee
By Theorem B it implies that $u=0$.\smallskip

\noindent If $\frac{N}{N-1}<p<\frac{N-1}{N-2}$, $\CE_+$ is discrete. Then either $v(t,.)\to \gw_s$ of $v(t,.)\to 0$ when $t\to-\infty$. In the first case 
it is equivalent to 
\bel{II-C-9}\BA {lll}
|x|^{\frac{1}{p-1}}u(x)=\gw_s(\frac{x}{|x|}) (1+o(1))\quad\text{as }x\to 0,
\EA\ee
and in the second case 
\bel{II-C-10}\BA {lll}\displaystyle
\lim_{x\to 0}|x|^{\frac{1}{p-1}}u(x)=0.
\EA\ee
\subsubsection{Weak singularities}
In the sequel, we assume $N>2$, the proof in the case $N=2$ can be carried out by the same techniques with minor technical modifications. 
\bprop{decay}If $(\ref{II-C-10})$ holds we claim that there exists $\gd>0$ such that 
\bel{II-C-11}\BA {lll}\displaystyle
u(x)\leq c|x|^{\gd-\frac{1}{p-1}},
\EA\ee
near $x=0$ for some $c>0$. 
\es

We proceed by contradiction, set $\gr(t)=\norm{v(t,.)}_{C^0{(S^{N-1}_+)}}$ and assume that for any $\ge>0$,
\bel{II-C-12}\BA {lll}\displaystyle
\limsup_{t\to-\infty}e^{-\ge t}\gr(t)=\infty.
\EA\ee

The following lemma proved in \cite{CMV} is the key for starting the proof of the decay of the solution.
\blemma{CMV)lem} There exists a function $\eta\in C^{\infty}((-\infty,T_1])$ satisfying
\bel{II-C-13}\BA {lll}\displaystyle
(i)\quad\quad\displaystyle&\eta>0\,,\;\eta_t>0\,,\;\displaystyle\lim_{t\to-\infty}\eta(t)=0,\\\displaystyle
(ii)\quad\quad &\displaystyle0<\limsup_{t\to-\infty}\myfrac{\gr(t)}{\eta(t)}<\infty,\\[3mm]
(iii)\quad\quad &\displaystyle\lim_{t\to-\infty}e^{-\ge t}\eta(t)=\infty\quad\text{for all }\ge>0,\\\displaystyle
(iv)\quad\quad &\displaystyle\myfrac{\eta_t}{\eta}\text{ and }\left(\myfrac{\eta_t}{\eta}\right)_t\text{ are bounded and integrable on  }(-\infty, T_1],\\\displaystyle
(v)\quad\quad &\displaystyle\lim_{t\to-\infty}\myfrac{\eta_t}{\eta}(t)=\lim_{t\to-\infty}\left(\myfrac{\eta_t}{\eta}\right)_t\!\!\!(t)=0.
\EA\ee
\es

\noindent{\it Proof of \rprop{decay}}. Define $w(t,.)=\eta^{-1}(t)v(t,.)$. Then $w$ is bounded and satisfies 
\bel{II-C-14}\BA{lll}
(1+\ge_1)w_{tt}+\left(N-\myfrac{2p}{p-1}+\ge_2+2(1+\ge_1)\myfrac{\eta_t}{\eta}\right)w_t+\Gd'w\\[2mm]\phantom{----}
+\left(\ell_{N,p}+\ge_3+(1+\ge_1)\myfrac{\eta_{tt}}{\eta}+\left(N-\myfrac{2p}{p-1}+\ge_2\right)\myfrac{\eta_t}{\eta}\right)w\\[3mm]
\phantom{----}+
\langle\nabla'w,{\ge_4}+\myfrac{\eta_t}{\eta}\ge_5\rangle+\langle\nabla'w_t,{\ge_5}\rangle
+\langle\nabla'\langle\nabla'w,{\bf e}_{_N}\rangle,{\ge_6}\rangle=0
\EA
\ee
in $(-\infty,T_0]\ti S^{N-1}_+$ and
\bel{II-C-14'}\BA{lll}
-\left(\left(w_t+\left(\myfrac{\eta_t}{\eta}-\myfrac{1}{p-1}\right)w\right)\langle{\bf c},{\bf e}_{_N}\rangle+\langle\nabla' w,{\bf e}_{_N}\rangle\right)\ge_7\\[3mm]
\phantom{----}
+\left(w_t+\left(\myfrac{\eta_t}{\eta}-\myfrac{1}{p-1}\right)w\right)\ge_8+\langle\nabla'w,{\gn}\rangle+\eta^{p-1}w^p=0
\EA
\ee
in $(-\infty,T_0]\ti \prt S^{N-1}_+$.
Since $w$ is bounded, a standard adaptation of \rlemma{est-sup} shows that $w_t,w_{tt}, w_{ttt}$ and $\nabla^\ga\prt_{t^\gb}w$ are uniformly bounded on 
$(-\infty,T_1]\ti S^{N-1}_+$ whenever $|\ga|+\gb\leq 3$ and $T_1<T_0$. The negative trajectory of $t\mapsto w(t,.)$ in $C^{2}(S^{N-1}_+)$ is defined by $\CT_-[w]:=\bigcup_{t\leq T_1}\{w(t,.)\}$. By the previous statements, $\CT_-[w]$ is bounded in $C^{3}(S^{N-1}_+)$, hence it is relatively compact in $C^{2}(S^{N-1}_+)$ by the Arzela-Ascoli theorem. Therefore, the alpha-limit set of $\CT_-[w]$ defined by 
\bel{II-C-15}
\CA[\CT_-[w]]:=\bigcap_{t\leq T_1}\text{clo}_{C^2(S^{N-1}_+)}\left(\bigcup_{\gt\leq t}\{w(\gt,.)\}\right)
\ee
is a non-empty compact connected set in $C^2(S^{N-1}_+)$. The integrability assumptions on $\eta$ allows us to prove
\blemma{energy-w} There holds
\bel{II-C-16}\BA {lll}
\myint{-\infty}{T_1}\myint{S^{N-1}_+}{}\left(w_t^2+w_{tt}^2+|\nabla'w_t|^2\right) dSdt<\infty.
\EA\ee
\es
\Proof We multiply equation $(\ref{II-C-14})$ by $w_t$ and integrate over $S^{N-1}_+$.
\bel{II-C-17}\BA {lll}
\myfrac{1}{2}\myfrac{d}{dt}\left[\myint{S^{N-1}_+}{}\left(w_t^2+\ell_{N,p}w^2-|\nabla'w|^2\right)dS-\myfrac{2\eta^{p-1}}{p+1}
\myint{\prt S^{N-1}_+}{}|w|^{p+1}dS'\right]\\[4mm]\phantom{--}
+\left(N-\myfrac{2p}{p-1}+\ge_2+2(1+\ge_1)\myfrac{\eta_t}{\eta}\right)\myint{S^{N-1}_+}{}w_t^2dS+\ga_1(t)+\ga_2(t)=0,
\EA\ee
where $\ga_1$ and $\ga_2$ are defined  by

$$\BA {lll}\ga_1(t)=\myint{S^{N-1}_+}{}\left[\ge_1w_{tt}+\ge_2w_t+\left(\ge_3+(1+\ge_1)\myfrac{\eta_{tt}}{\eta}+\left(N-\myfrac{2p}{p-1}+\ge_2\right)\myfrac{\eta_t}{\eta}\right)w\right.\\[4mm]\phantom{-------}
\left.+\langle\nabla'w,{\ge_4}+\myfrac{\eta_t}{\eta}\ge_5\rangle+\langle\nabla'w_t,\ge_5\rangle+\langle\nabla'\langle\nabla'w,{\bf e}_{_N}\rangle,\ge_6\rangle\right]w_tdS,
\EA$$
$$\BA{lll}\ga_2(t)=\myint{\prt S^{N-1}_+}{}\left[\left(w_t+\left(\myfrac{\eta_t}{\eta}-\myfrac{1}{p-1}\right)w\right)\left(\langle{\bf c},{\bf e}_{_N}\rangle)\ge_7+\ge_8\right)\right.\\[4mm]\phantom{---------------}\left.
-\langle\nabla' w,{\bf e}_{_N}\rangle\ge_7-\myfrac{p-1}{p+1}w^{p+1}\eta^{p-1}\myfrac{\eta_t}{\eta}\right] w_tdS'.
\EA$$
 Using the estimates on $\ge_j$ and $(\ref{II-C-13})$, we obtain that
 $$\myint{-\infty}{T_1}\myint{S^{N-1}_+}{}w_t^2 dSdt<\infty.
 $$
 Multiplying equation $(\ref{II-C-14})$ by $w_{tt}$, integrating over $S^{N-1}_+$ and using $(\ref{II-C-13})$ yield 
$$\myint{-\infty}{T_1}\myint{S^{N-1}_+}{}\left(w_{tt}^2+|\nabla'w_t|^2\right) dSdt<\infty,$$
which ends the proof.\qeda\medskip

\noindent{\it End of the proof of \rprop{decay}}. Since $w_t$ and $w_{tt}$ are uniformly continuous on $(-\infty,T_1]\ti S^{N-1}_+$ we infer from $(\ref{II-C-16})$ that 
\bel{II-C-18}\BA {lll}
\displaystyle \lim_{t\to-\infty}w_t(t,.)=\lim_{t\to-\infty}w_{tt}(t,.)=0,
\EA\ee
uniformly on $S^{N-1}_+$. Therefore $\CA[\CT_-[w]]$ is a subset of the set of nonnegative solutions of 
\bel{II-C-19}\BA {lll}
\Gd'\gf+\ell_{N,p}\gf=0\qquad&\text{in }S^{N-1}_+\\[2mm]
\phantom{-+,}
\langle\nabla'\gf,{ \gn}\rangle=0
\qquad&\text{in }\prt S^{N-1}_+,
\EA\ee
and by $(\ref{II-C-13})$-(ii) it contains a positive element. Since $\ell_{N,p}>0$ this is a contradiction and $(\ref{II-C-12})$ does not hold. This ends the proof of
\rprop{decay}.
 $\phantom{----------}$\qeda\medskip

\noindent{\it Step 1}. We claim that 
\bel{II-C-19'}
u(x)\leq c|x|^{2-N}\quad\text {in a neighborhood of } 0.
\ee
If $\gd\geq\frac{1}{p-1}+2-N$, $(\ref{II-C-19'})$ is a consequence of $(\ref{II-C-11})$. In what follows we assume that
\bel{II-C-19"}
0<\gd<\frac{1}{p-1}+2-N.
\ee
We set $v_{\gd}=e^{-\gd t}v$. Then $v_{\gd}$ is bounded in $(-\infty,T_1]\ti S^{N-1}_+$ and, as in the proof of \rprop{decay}, the quantities $\prt_tv_{\gd}$, $\prt_{tt}v_{\gd}$, $\prt_{ttt}v_{\gd}$ and $\nabla^\ga\prt_{t^\gb}v_{\gd}$ are uniformly bounded on $(-\infty,T_1]\ti S^{N-1}_+$. Furthermore there holds
\bel{II-C-20}\BA{lll}
(1+\ge_1)v_{\gd\,tt}+\left(N-\myfrac{2p}{p-1}+\ge_2+2(1+\ge_1)\gd\right)v_{\gd\,t}+\Gd'v_{\gd}\\[2mm]\phantom{----}
+\left(\ell_{N,p}+\ge_3+(1+\ge_1)\gd^2+\left(N-\myfrac{2p}{p-1}+\ge_2\right)\gd\right)v_{\gd}\\[3mm]
\phantom{----}+
\langle\nabla'v_{\gd},{\ge_4}+\gd\ge_5\rangle+\langle\nabla'v_{\gd\,t},{\ge_5}\rangle
+\langle\nabla'\langle\nabla'v_{\gd},{\bf e}_{_N}\rangle,{\ge_6}\rangle=0
\EA
\ee
in $(-\infty,T_0]\ti S^{N-1}_+$ and
\bel{II-C-21}\BA{lll}
-\left(\left(v_{\gd\,t}+\left(\gd-\myfrac{1}{p-1}\right)v_{\gd}\right)\langle{\bf c},{\bf e}_{_N}\rangle+\langle\nabla' v_{\gd},{\bf e}_{_N}\rangle\right)\ge_7\\[3mm]
\phantom{----}
+\left(v_{\gd\,t}+\left(\gd-\myfrac{1}{p-1}\right)v_{\gd}\right)\ge_8+\langle\nabla'v_{\gd},{\gn}\rangle+e^{(p-1)t}v_{\gd}^p=0
\EA
\ee
in $(-\infty,T_0]\ti \prt S^{N-1}_+$. We denote by $\tilde v_\gd$ the projection of $v$ onto $H:=[\ker(-\Gd')]^\perp$, the operator being defined in $H^1(S^{N-1}_+)$, and by $P_H$ the corresponding projection operator. Then 
  \bel{II-C-22}\BA {lll}
\tilde v_{\gd\,tt}+\left(N-\myfrac{2p}{p-1}+2\gd\right)\tilde v_{\gd\,t}+\Gd'\tilde v_{\gd}\\[2mm]\phantom{----}
+\left(\ell_{N,p}+\gd^2+\left(N-\myfrac{2p}{p-1}\right)\gd\right)\tilde v_{\gd}+\tilde F=0
\EA\ee
where 
$$\BA{lll}\Tilde F:=P\left[\ge_1v_{\gd\,tt}+(\ge_2+2\ge_1\gd)v_{\gd\,t}+(\ge_3+\ge_1\gd^2+\ge_2\gd)v_\gd\right.\\[2mm]
\phantom{------}\left.+\langle\nabla'v,{\ge_4}\rangle+\langle\nabla'v_t,{\ge_5}\rangle
+\langle\nabla'\langle\nabla'v,{\bf e}_{_N}\rangle,{\ge_6}\rangle\right]=O(e^t).
\EA$$
We multiply $(\ref{II-C-21})$ by $\tilde v_\gd$, integrate over $S^{N-1}_+$ and use the boundary condition  and  the fact that $N-1$ is the first eigenvalue of $-\Gd'$ in $H$. We deduce that 
$\tilde X_\gd(t):=\norm{\tilde v_\gd(t,.)}_{L^2(S^{N-1}_+)}$ satisfies in the sense of distributions in $(-\infty,T_1)$,
  \bel{II-C-23}\BA {lll}
\tilde X''_\gd+\left(N-\myfrac{2p}{p-1}+2\gd\right)\tilde X'_\gd\\[2mm]
\phantom{-----}
+\left(\ell_{N,p}+\gd^2+\left(N-\myfrac{2p}{p-1}\right)\gd+1-N\right)\tilde X_\gd\geq -c^*e^{mt}.
\EA\ee
where $m=\inf\{1,p-1\}$, for some constant $c^*>0$. Note that the nonlinear term on $\prt S^{N-1}_+$ is at the origin of the term $e^{(p-1)t}$. The characteristic polynomial of $(\ref{II-C-23})$ is
$$P_\gd (\xi)=\xi^2+\left(N-\myfrac{2p}{p-1}+2\gd\right)\xi+\ell_{N,p}+\gd^2+\left(N-\myfrac{2p}{p-1}\right)\gd+1-N.
$$
It is noticeable that its discriminant  is $N^2$, independent of $\gd$,  and as a consequence its roots are expressed easily by
  \bel{II-C-25'}\BA {lll}
\xi_{1,\gd}=\myfrac{p}{p-1}-\gd>m\quad\text{and }\; \xi_{2,\gd}=\myfrac{p}{p-1}-N-\gd<0,
\EA\ee
since $(\ref{II-C-19"})$ holds. Therefore $P_\gd(m)<0$.  For $a,\gg,\ge>0$ set  
$$X_\ge(t)= ae^{\xi_{1,\gd}t}+\ge e^{\xi_{1,\gd}t} +\gg e^{mt}.
$$
Then 
$$\BA {lll}
X''+\left(N-\myfrac{2p}{p-1}+2\gd\right)X'\\
[2mm]\phantom{-----}
+\left(\ell_{N,p}+\gd^2+\left(N-\myfrac{2p}{p-1}\right)\gd+1-N\right)X=\gg P_\gd(m)e^{mt}.
\EA$$
We can choose $\gg$ such that $\gg P_\gd(1)\geq -c^*$ and $a=\norm{\tilde v_\gd(T_1,.)}_{L^2(S^{N-1}_+)}e^{-\xi_{1,\gd}T_1}$. 
By the maximum principle $\tilde X_\gd(t)\leq X_\ge(t)$ for $t\leq T_1$ and all $\ge>0$. This implies 
\bel{II-C-26}\BA {lll}
\norm{\tilde v_\gd(t,.)}_{L^2(S^{N-1}_+)}\leq \norm{\tilde v_\gd(T_1,.)}_{L^2(S^{N-1}_+)}e^{\xi_{1,\gd}(t-T_1)}+\gg e^{mt}
\quad\text{for }t\leq T_1.
\EA\ee
Using standard regularizing effect for elliptic equations, we can improve $(\ref{II-C-26})$ and obtain a uniform estimate
\bel{II-C-27}\BA {lll}
\norm{\tilde v_\gd(t,.)}_{L^\infty(S^{N-1}_+)}\leq A\norm{\tilde v_\gd(T_1,.)}_{L^2(S^{N-1}_+)}e^{\xi_{1,\gd}(t-T_1)}+\gg e^{mt}
\quad\text{for }t\leq T_1-1.
\EA\ee

Next we denote by $X_\gd$ the projection of $v_\gd$ onto $\ker(-\Gd')$ (i.e. the average on $S^{N-1}_+$), then 
 \bel{II-C-28}\BA {lll}
X''_\gd+\left(N-\myfrac{2p}{p-1}+2\gd\right)X'_\gd
+\left(\ell_{N,p}+\gd^2+\left(N-\myfrac{2p}{p-1}\right)\gd\right)X_{\gd}+ F=0
\EA\ee
where 
$$\BA{lll} F:=\myfrac{1}{|S^{N-1}_+|}\myint{S^{N-1}_+}{}\left[\ge_1v_{\gd\,tt}+(\ge_2+2\ge_1\gd)v_{\gd\,t}+(\ge_3+\ge_1\gd^2+\ge_2\gd)v_\gd\right.\\[2mm]
\phantom{------}\left.+\langle\nabla'v,{\ge_4}\rangle+\langle\nabla'v_t,{\ge_5}\rangle
+\langle\nabla'\langle\nabla'v,{\bf e}_{_N}\rangle,{\ge_6}\rangle\right]dS=O(e^{mt}).
\EA$$
The characteristic roots of the equation
$$y''+\left(N-\myfrac{2p}{p-1}+2\gd\right)y'
+\left(\ell_{N,p}+\gd^2+\left(N-\myfrac{2p}{p-1}\right)\gd\right)y=0
$$
are $\gth_{1,\gd}$, $\gth_{2,\gd}$. They can easily be computed and for $\gd>0$ small enough
 \bel{II-C-31}\gth_{1,\gd}=\myfrac{1}{p-1}-\gd>N-2\geq 1>\gth_{2,\gd}=\myfrac{1}{p-1}+2-N-\gd.
\ee 
The solution of $(\ref{II-C-28})$ admits the general expression
 \bel{II-C-32}
X_{\gd}(t)=ae^{t\gth_{1,\gd}}+be^{t\gth_{2,\gd}}-
\myfrac{1}{\gth_{1,\gd}-\gth_{2,\gd}}\myint{t}{T_1}F(s)\left(e^{(t-s)\gth_{1,\gd}}-e^{(t-s)\gth_{2,\gd}}\right)ds
\ee
Since $m<\gth_{1,\gd}$, it is easy to see that there exists $c\geq 0$ such that, when $t\to-\infty$, there holds
 \bel{II-C-33}\BA {lll}
X_\gd(t)=e^{t\inf\{p-1,\gth_{2,\gd}\}}(c+o(1)),
\EA\ee
if $p-1\neq \gth_{2,\gd}$ and 
 \bel{II-C-34}\BA {lll}
X_\gd(t)=(-t)e^{t(p-1)}(c+o(1)),
\EA\ee
if $p-1= \gth_{2,\gd}$. We consider only the case $p-1\neq \gth_{2,\gd}$, the case of equality requiring only some technical modifications of the proof. \smallskip

\noindent {\it Case 1:} Assume $\gth_{2,\gd}<p-1$. Then $X_\gd(t)=e^{t\gth_{2,\gd}\}}(c+o(1))$. Since $\xi_{1,\gd}>1$ and $m>\gth_{2,\gd}$, 
we infer from $(\ref{II-C-27})$ and $(\ref{II-C-31})$ and the definition of $v_\gd$ that 
 \bel{II-C-35}\BA {lll}
0<v(t,\gs)= e^{t(\gd+\gth_{2,\gd})}(c+o(1))\quad\text{when }t\to-\infty,\text{ uniformly on }S^{N-1}_+.
\EA\ee
for some constant $c>0$. This implies not only $(\ref{II-C-19'})$ but also $(\ref{I-1-13})$.\smallskip

\noindent {\it Case 2:} Assume $1>\gth_{2,\gd}>p-1=m$. Then 
 \bel{II-C-36}\BA {lll}
0<v(t,\gs)\leq c e^{(\gd+p-1)t}\quad\text{for all }t\in (-\infty,T_1]\ti S^{N-1}_+.
\EA\ee
Then we restart the previous construction, replacing $\gd$ by $\gd_1:=\gd+p-1$. After a finite number $j$ of iterations of this construction 
and setting $\gd_j:=\gd+j(p-1)$ we finally obtain
 \bel{II-C-37}\BA {lll}
0<v(t,\gs)= e^{t(\gd_j+\gth_{2,\gd_j})}(c+o(1))\quad\text{when }t\to-\infty,\text{ uniformly on }S^{N-1}_+.
\EA\ee
which again implies not only $(\ref{II-C-19'})$ but also $(\ref{I-1-13})$.\qeda\medskip

\noindent \Remark The results of Theorem C can be extended to signed solutions $u$ of $(\ref{I-1-9})$ provided they satisfy not only the same a priori estimates
$|u(x)|\leq c|x|^{-\frac{1}{p-1}}$ but also  $|D^\ga u(x)|\leq c|x|^{-\frac{1}{p-1}-|\ga|}$ for $|\ga|=1,2,3$. If this holds, the energy method applies and we infer that 
the limit set of the trajectory $\CA[\CT_-[v]]$ is a connected subset of the set $\CE$. In particular, if $p\geq \frac{N-1}{N-2}$ then $\CA[\CT_-[v]]=\{0\}$ and by Theorem B it implies that $u=0$. If $\frac{N}{N-1}<p<\frac{N-1}{N-2}$ with $N>2$, then $\CA[\CT_-[v]]\subset\{\gw_s,-\gw_s,0\}$. If $N=2$ and 
$\frac{1}{p-1}$ is an integer, then $\CA[\CT_-[v]]=\{0\}$ and  $u=0$ by Theorem B, while if $\frac{1}{p-1}$ is not an integer then $\CA[\CT_-[v]]\subset\{\gw_s,-\gw_s,0\}$. Furthermore, if $\CA[\CT_-[v]]=\{0\}$ and $\ell_{N,p}$ is not an eigenvalue of $-\Gd'$ in $H^1(S^{N-1}_+)$ it is possible to adapt the method developed in the proof of \rprop {decay} and obtain that $r^{N-2+k}u(r,.)$ converges to a nonzero eigenfunction of $-\Gd'$ in $H^1(S^{N-1}_+)$ for some 
$k\in\BBN$ such that $N-2+k<\frac{1}{p-1}$. The method for such a task is an adaptation of the ideas introduced in \cite[Theorem 2.1]{CMV} and \cite[Theorem 5.1]{GmVe}. Note that the assumption $\ell_{N,p}\notin\gs\left(-\Gd',H^1(S^{N-1}_+)\right)$ is fundamental to prove and estimate of type $(\ref{II-C-11})$, which is the starting point of the proof. 
\mysection{Measure boundary data}
Let $u$ and $v$ two solutions of $(\ref{I-1-0})$ with the same data $\gm$. By \rlemma{lem-br} 
\bel{III-1-1}
\myint{\Gw}{}|u-v|\left(-\Gd \gz+\gz\right) dx+\myint{\prt\Gw}{}{\rm{sign}}_0(u-v)(g(u)-g(v)\gz dS\leq 0
\ee
for all $\gz\in\CC(\Gw)$, $\gz\geq 0$. Since $g$ is nondecreasing, we take $\gz=1$ and get $u=v$. 
\subsection{Proof of Theorem D}
In this section we assume $N\geq 3$. Let $\{\gm_k\}$ be a sequence of smooth functions on $\prt\Gw$ and $u_k$ the solution of 
\bel{III-1-2}\BA {lll}
\phantom{g(u;}
-\Gd u_k+u_k=0\qquad&\text{in }\;\Gw\\[0mm]\phantom{,,,u}
\myfrac{\prt u_k}{\prt {\bf n}}+g(u_k)=\gm_k\qquad&\text{in }\;\prt\Gw,
\EA
\ee
obtained by minimization. By \rlemma{lem-br}
\bel{III-1-3}\BA {lll}
\myint{\Gw}{}|u_k|dx+\myint{\prt\Gw}{}|g(u_k)| dS\leq \myint{\prt\Gw}{}|\gm_k|dS.
\EA
\ee
Hence
\bel{III-1-4}\BA {lll}
\norm{u_k}_{L^{\frac{N}{N-2},\infty}(\Gw)}+\norm{\nabla u_k}_{L^{\frac{N}{N-1},\infty}(\Gw)}\leq c'\myint{\prt\Gw}{}|\gm_k|dS\leq c
\norm{\gm}_{\mathfrak M(\prt\Gw)},
\EA
\ee
by \rprop {lin}.

Therefore there exist a function $u\in L^{\frac{N}{N-2},\infty}(\Gw)$ verifying 
$\nabla u\in L^{\frac{N}{N-1},\infty}(\Gw)$  and a subsequence $\{u_{k_j}\}$ such that $u_{k_j}\to u$ a.e. in $\Gw$ and  in $L^1(\Gw)$. By 
\cite{Triebel} the boundary trace of a function $v\in L^{\frac{N}{N-1},\infty}(\Gw)$ such that $\nabla v\in L^{\frac{N}{N-1},\infty}(\Gw)$ belongs to the fractional Besov-Lorentz space $B^{\frac{1}{N},\frac{N}{N-1},\infty}(\prt\Gw)$ and there holds
\bel{III-1-5}\BA {lll}
\norm{v\lfloor_{\prt\Gw}}_{B^{\frac{1}{N},\frac{N}{N-1},\infty}(\prt\Gw)}\leq c\left(\norm{v}_{L^{\frac{N}{N-1},\infty}(\Gw)}+\norm{\nabla v}_{L^{\frac{N}{N-1},\infty}(\Gw)}\right). 
\EA
\ee
Using Sobolev imbedding theorem for Besov-Lorentz spaces, classicaly obtained by the real interpolation method \cite{LP} from the same indexed Sobolev spaces \cite{Triebel}, we obtain
\bel{III-1-6}\BA {lll}
\norm {v\lfloor_{\prt\Gw}}_{L^{\frac{N-1}{N-2},\infty}(\prt\Gw)}\leq c\norm{v\lfloor_{\prt\Gw}}_{B^{\frac{1}{N},\frac{N}{N-1},\infty}(\prt\Gw)}.
\EA
\ee
Therefore 
\bel{III-1-50}\BA {lll}
\norm {u_k\lfloor_{\prt\Gw}}_{L^{\frac{N-1}{N-2},\infty}(\prt\Gw)}\leq c
\norm{\gm}_{\mathfrak M(\prt\Gw)},
\EA
\ee
and $u_{k_j}\lfloor_{\prt\Gw}\to u\lfloor_{\prt\Gw}$ a.e. in $\prt\Gw$. In order to prove the convergence of $\{g(u_{k_j})\lfloor_{\prt\Gw}\}$ to 
$\{g(u)\lfloor_{\prt\Gw}\}$ we use Vitali's theorem. Let $E\subset\prt\Gw$ be a Borel set and $|E|_{_{N-1}}=\CH_{_{\prt\BBR^N_+}}(E)$ is its the (N-1)-Hausdorff measure on $\prt\Gw$; for any $\gl>0$, 
\bel{III-1-51}\BA {llll}
\myint{E}{}|g(u_{k_j})|dS=\myint{E\cap\{|u_{k_j}|\leq \gl\}}{}|g(u_{k_j})|dS+ \myint{E\cap\{|u_{k_j}|>\gl\}}{}|g(u_{k_j})|dS\\[4mm]
\phantom{\myint{E}{}|g(u_{k_j})|dS}
\leq |E|_{_{N-1}}\left(g(\gl)-g(-\gl)\right)+\myint{\{|u_{k_j}\lfloor_{\prt\Gw}|>\gl\}}{}|g(u_{k_j})|dS. 
\EA\ee
We set $A_\gl(u_{k_j})=\{x\in\prt\Gw:|u_{k_j}\lfloor_{\prt\Gw}(x)|>\gl\}$ and $\ga_{k_j}(\gl)=|A_\gl(u_{k_j})|_{_{N-1}}$. Since $(\ref{III-1-50})$ holds, 
$$\ga_{k_j}(\gl)\leq c
\norm{\gm}_{\mathfrak M(\prt\Gw)}\gl^{-\frac{N-1}{N-2}}. $$
Using Cavalieri's formula \cite{Cav},
\bel{III-1-52}\BA {lll}
\myint{\{|u_{k_j}|>\gl\}}{}|g(u_{k_j})|dS=-\myint{\gl}{\infty}g(s)d\ga_{k_j}(s)\\[4mm]\phantom{\myint{\{|u_{k_j}|>\gl\}}{}|g(u_{k_j})|dS}
\leq c\frac{N-1}{N-2}
\norm{\gm}_{\mathfrak M(\prt\Gw)}\myint{\gl}{\infty}\left(g(s)-g(-s)\right)s^{-\frac{2N-3}{N-2}}ds.
\EA\ee
Combining $(\ref{III-1-51})$ and $(\ref{III-1-52})$, we can choose $\gl$ large enough and  deduce that $\myint{E}{}|g(u_{k_j})|dS\to 0$ when $|E|_{_{N-1}}\to 0$, uniformly with respect to $k_j$. Hence $g(u_{k_j})\lfloor_{\prt\Gw}\to g(u)\lfloor_{\prt\Gw}$ in $L^1(\prt\Gw)$. If $\gx\in \CC(\Gw)$, there holds
\bel{III-1-10}\BA {lll}
\myint{\Gw}{}u_{k_j}\left(-\Gd\xi+\xi\right) dx+\myint{\prt\Gw}{}g(u_{k_j})\xi dS=\myint{\prt\Gw}{}\xi \gm_{k_j}dS.
\EA\ee
Letting $k_j\to\infty$, we infer that $(\ref{I-1-16})$ holds. Actually, the whole sequence $\{u_k\}$ converges and we denote  by $u_\gm$ its limit. 
Notice also that by the monotonicity of $g$, $\gm\geq\gm'$ implies $u_{\gm}\geq u_{\gm'}$.\qeda\medskip

\noindent \Remark If $g(r)=|r|^{p-1}r$ with $p>0$, condition $(\ref{I-1-20})$ is satisfied if and only if $p<\frac{N-1}{N-2}$. 
\subsection{Proof of Theorem E}
In this section we assume $N=2$. \medskip

\noindent {\it Proof of assertion 1}. As in the proof of Theorem D, we denote by 
$u_k$ the solution of $(\ref{III-1-2})$. Estimate $(\ref{III-1-3})$ is valid and $(\ref{III-1-4})$ is replaced by 
\bel{III-1-4}\BA {lll}
\norm{u_k}_{L^{q}(\Gw)}+\norm{\nabla u_k}_{L^{2,\infty}(\Gw)}\leq c'_q\myint{\prt\Gw}{}|\gm_k|dS\leq c_q\norm{\gm}_{\mathfrak M(\prt\Gw)}
\quad\text{for all }q\in [1,\infty).
\EA
\ee
By an extension of Moser's inequality to Lorentz spaces  \cite[Theorem 3.1]{XZ} , there exist  constants $c^*,c'>0$ depending on $\Gw$ such that for any function $v\in L^2(\prt\Gw)$ such that $(-\Gd)^{\frac{1}{2}}v\in L^{2,\infty}(\prt\Gw)$ (equivalently $v\in B^{\frac12,2,\infty}(\prt\Gw)$), there holds
\bel{III-1-5}\BA {lll}\displaystyle
\sup_{\norm{(-\Gd)^{\frac12}v}_{L^{2,\infty}}\leq 1}\myint{\prt\Gw}{}e^{c^* |v(x)|}dS\leq c'.
\EA
\ee 
Using $(\ref{III-1-4})$, $(\ref{III-1-5})$ we deduce, with $c=\frac{c^*}{c_2}$,
\bel{III-1-6}\BA {lll}\displaystyle
\myint{\prt\Gw}{}e^{\frac{c|u_k(x)|}{\norm{\gm}_{\mathfrak M}} }dS\leq c'.
\EA
\ee
This implies that $\{(u_k,u_k\lfloor_{\prt\Gw})\}$ is compact in $L^q(\Gw)\ti L^q(\prt\Gw)$ for  any $q<\infty$ and up to a subsequence $\{u_{k_j}\}$  converges a. e. and in $L^q(\Gw)\ti L^q(\prt\Gw)$ to some $u$ such that $\nabla u\in L^{2,\infty}(\Gw)$ and therefore $u\lfloor_{\prt\Gw}\in L^q(\prt\Gw)$. Thus $u\lfloor_{\prt\Gw}$ satisfies   $(\ref{III-1-6})$. As a consequence problem $(\ref{I-1-0})$ admits a solution if $|g(r)|\leq c_1|r|^q+c_2$ for some $q\in (0,\infty)$ and $c_1,c_2\geq 0$. We have actually a more general result if we assume that $a_+(g)=a_-(g)=0$. From $(\ref{III-1-6})$ there holds for $\gl>0$,
\bel{III-1-7}e^{\frac{c\gl}{\norm{\gm}_{\mathfrak M}} }\ga_{k_j}(\gl)\leq c'\Longrightarrow \ga_{k_j}\leq c'e^{-\frac{c\gl}{\norm{\gm}_{\mathfrak M}} },
\ee
where $A_\gl(u_{k_j})$ and $\ga_{k_j}(\gl)$ are defined in the proof of Theorem D. If $E\subset\prt\Gw$ is a Borel set, 
\bel{III-1-7}\BA {lll}\myint{E}{}|g(u_{k_j})|dS\leq |E|_1(g(\gl)-g(-\gl))+\myint{\{|u_{k_j}|\lfloor_{\prt\Gw}\}}{}|g(u_{k_j})|dS\\[4mm]
\phantom{\myint{E}{}|g(u_{k_j})|dS}
\leq |E|_1(g(\gl)-g(-\gl))-\myint{\gl}{\infty}(g(s)-g(-s))d\ga_{k_j}(s)\\[4mm]
\phantom{\myint{E}{}|g(u_{k_j})|dS}
\leq  |E|_1(g(\gl)-g(-\gl))+\myfrac{c'}{\norm{\gm}_{\mathfrak M}}\myint{\gl}{\infty}(g(s)-g(-s))e^{-\frac{cs}{\norm{\gm}_{\mathfrak M}} }ds
\EA\ee
Since 
$$\myint{\gl}{\infty}(g(s)-g(-s))e^{-as} ds<\infty$$
for any $a>0$ the result follows as in Theorem D. \qeda\medskip

\noindent \Remark Actually we have a stronger result since we only use
\bel{III-1-8}
\myint{\gl}{\infty}(g(s)-g(-s))e^{-\frac{cs}{\norm{\gm}_{\mathfrak M}} }ds<\infty.
\ee
Therefore the assumption $a_+(g)=a_-(g)=0$ can be replaced by $\norm{\gm}_{\mathfrak M}\leq c_g$ where 
\bel{III-1-9}
\myint{\gl}{\infty}(g(s)-g(-s))e^{-\frac{cs}{c_g} }ds<\infty.
\ee
However the constant $c$ is not explicitely known.\medskip

\noindent {\it Proof of assertion 2}. Set $\displaystyle\gm=\sum_{j=1}^k\ga_j\gd_{a_j}$. For $\ell>0$ set $g_\ell(r)=\min\{g(\ell),\sup\{g(-\ell),g(r)\}\}$. Since $a_+(g_\ell)=a_-(g_\ell)=0$, there exists a weak solution  to 
\bel{III-1-10}\BA {lll}
\phantom{,}-\Gd u+u=0\qquad&\text{in }\,\Gw\\[0mm]
\myfrac{\prt u}{\prt{\bf n}}+g_\ell(u)=\gm&\text{on }\,\prt\Gw,
\EA\ee
and this solution denoted by $u_{\ell,\gm}$ is unique since $g_\ell$ is nonnecreasing. Put $J_+:=\{j=1,...,k:\ga_j>0\}$, $J_-:=\{j=1,...,k:\ga_j<0\}$ and denote by 
$u_{\ell,\gm_+}$ (resp. $u_{\ell,\gm_-}$) the solution of 
\bel{III-1-10'}\BA {lll}
\phantom{\myfrac{\prt u}{\prt{\bf n}}+}-\Gd u + u =0\qquad&\text{in }\,\Gw\\[0mm]
\myfrac{\prt u}{\prt{\bf n}}+g_\ell(u)=\gm_+\quad\text{(resp. $=\gm_-$)}&\text{on }\,\prt\Gw,
\EA\ee
with $\displaystyle\gm_+\!=\!\!\sum_{j\in J_+}\ga_j\gd_{a_j}$ (resp. $\displaystyle\gm_-\!=\!\!\sum_{j\in J_-}\ga_j\gd_{a_j}$). Then $u_{\ell,\gm_+}\geq 0$ (resp. $u_{\ell,\gm_-}\leq 0$) and
$$\sum_{j\in J_-}\ga_j\CN_\Gw(.,a_j)\leq u_{\ell,\gm_-}\leq u_{\ell,\gm}\leq u_{\ell,\gm_+}\leq \sum_{j\in J_+}\ga_j\CN_\Gw(.,a_j).
$$
Thus
$$\BA {lll}\displaystyle g\left(\sum_{j\in J_-}\ga_j\CN_\Gw(.,a_j)\right)\leq g_{\ell}\left(u_{\ell,\gm_-}\right) \leq g_\ell(u_{\ell,\gm})\leq g_{\ell}\left(u_{\ell,\gm_+}\right)\leq g\left(\sum_{j\in J_+}\ga_j\CN_\Gw(.,a_j)\right).
\EA$$
Since
$$\ga_j\CN_\Gw(x,a_j)=\myfrac{\ga_j}{\gp}\ln\left(\myfrac{1}{|x-a_j|}\right)(1+o(1))\text {as }\,x\to a_j,
$$
for any $\ge>0$, there exists $K_\ge>0$ such that 
$$\BA {lll}\displaystyle\sum_{j\in J_-}g\left(\myfrac{\ga_j-\ge}{\gp}\ln\left(\myfrac{1}{|x-a_j|}\right)\right)-K_\ge\leq g_{\ell}\left(u_{\ell,\gm_-}\right) \leq g_\ell(u_{\ell,\gm})\\[4mm]\phantom{------}\displaystyle\leq g_{\ell}\left(u_{\ell,\gm_+}\right)\displaystyle\leq \sum_{j\in J_+}g\left(\myfrac{\ga_j+\ge}{\gp}\ln\left(\myfrac{1}{|x-a_j|}\right)\right)+K_\ge.
\EA$$
We take $\ge>0$ small enough such that 
$$\displaystyle\myfrac{\gp}{a_-(g)}\inf_{j\in J_-}\ga_j-\ge\leq 0\leq\sup_{j\in J_+}\ga_j+\ge\leq \myfrac{\gp}{a_+(g)}.
$$
This implies that $\{g_{\ell}\left(u_{\ell,\gm_-}\right)\}_\ell$ and $\{g_{\ell}\left(u_{\ell,\gm_+}\right)\}_\ell$ are uniformly integrable in $L^1(\prt\Gw)$.  Consequently 
$\{g_{\ell}\left(u_{\ell,\gm}\right)\}_\ell$ is also uniformly integrable in $L^1(\prt\Gw)$. Letting $\ell\to \infty$ we deduce that up to a subsequence, $u_{\ell_j,\gm}$ converges to the unique weak solution $u=u_\gm$ of $(\ref{I-1-9})$.\qeda\medskip

\noindent \Remark By adapting the construction in \cite{Vazquez} (see also \cite{Veron} for a slightly simpler proof), it can be proved that when $N=2$
the problem  $(\ref{I-1-9})$  can be solved with any measure on $\prt\Gw$ with Jordan decomposition $\gm=\gm_r+\gm_a$ where $\gm_r$ is the diffuse part and $\gm_a=\sum_{j=1}^k\ga_j\gd_{a_j}$ is the atomic part, provided  the $\ga_j$ satisfy $(\ref{-1-22})$. In particular no assumption on $\gm_r$ is required. 
\subsection{The supercritical case: proof of Theorem F}
Let  $\BBP_\Gw$ be the Poisson operator for $-\Gd+I$ in $\Gw$ and $\BBD_\Gw$  the Dirichlet to Neumann operator for  $-\Gd+I$. Thus if $\eta\in\mathfrak M(\prt\Gw)$, $v=\BBP_\Gw[\eta]$ if 
\bel{IV-1-1}\BA {lll}
-\Gd v+v=0\qquad&\text{in }\,\Gw\\[0mm]
\phantom{-\Gd +v}
v=\eta&\text{on }\,\prt\Gw,
\EA\ee
and 
\bel{IV-1-2}\BA {lll}
\BBD_\Gw[\eta]=\myfrac{\prt v}{\prt{\bf n}}=\myfrac{\prt }{\prt{\bf n}}\BBP_\Gw[\eta]\quad\text{on }\,\prt\Gw.
\EA\ee
Let $\BBN_\Gw$ be the Neumann operator from $\prt\Gw$ to $\Gw$ defined by $v=\BBN_\Gw[\gm]$ where
\bel{IV-1-0}\BA {lll}
-\Gd v+v=0\qquad&\text{in }\,\Gw\\[0mm]
\phantom{-\Gd-,}
\myfrac{\prt v}{\prt {\bf n}}=\gm&\text{on }\,\prt\Gw,
\EA\ee
and some results of regularity of $\BBN_\Gw$ are recalled in \rprop{QR}.
\blemma{regDN} Let $1<q<\infty$ and $\gm$ is a distribution in $\BBR^N$ with support included in $\prt\Gw$. Then the following assertions are equivalent:
\smallskip

\noindent (i) $\BBN_\Gw[\gm]\in L^q(\prt\Gw)$,\smallskip

\noindent (ii) $\gm\in W^{-1,q}(\prt\Gw)$. \smallskip

\noindent Furthermore there exists $c>0$ such that 
\bel{IV-1-3}\BA {lll}
c^{-1}\norm\gm_{W^{-1,q}(\prt\Gw)}\leq \norm{\BBN_\Gw[\gm]}_{L^q(\prt\Gw)}\leq c\norm\gm_{W^{-1,q}(\prt\Gw)}
\EA\ee
\es
\Proof We recall that by Calderon's theorem the operator $\BBD_\Gw$ is an isomorphism from $L^q(\prt\Gw)$ to $W^{1,q}(\prt\Gw)$ (see e.g. \cite[Theorem 1.2.3]{AH}) and in particular for any 
$q\in (1,\infty)$, 
\bel{IV-1-4}\BA {lll}
c^{-1}\norm\eta_{W^{1,q}(\prt\Gw)}\leq \norm{\BBD_\Gw[\eta]}_{L^q(\prt\Gw)}\leq c\norm\eta_{W^{1,q}(\prt\Gw)}.
\EA\ee
This follows from the fact the following identity holds  if $v=\BBN_\Gw[\gm]$
$$v(x)=\myint{\prt\Gw}{}\CN_\Gw(x,y)\myfrac{\prt v}{\prt {\bf n}}(y)dS(y)=\myint{\prt\Gw}{}\CN_\Gw(x,y)d\gm(y)
$$
and that $\CN_\Gw$ satisfies $(\ref{I-B-1})$. In the flat case it is exactely the Calderon-Zygmund theory as it is shown in \cite[Theorem V-3]{Stein}.\smallskip

\noindent Let $v=\BBN_\Gw[\gm]$,  $\xi\in C^2(\Gw)$ and $\psi=\BBP_\Gw[\xi]$, then 
$$\langle \gm,\xi\rangle=\myint{\prt\Gw}{}v\lfloor_{\prt\Gw}\myfrac{\prt \psi}{\prt {\bf n}}dS=
\myint{\prt\Gw}{}\BBN_\Gw[\gm]\BBD_\Gw[\xi]dS.
$$
Using (\ref{IV-1-4}), we see that if $\BBN_\Gw[\gm]$ belongs to $L^q(\prt\Gw)$, then 
$$|\langle \gm,\xi\rangle|\leq \norm{\BBN_\Gw[\gm]}_{L^q(\prt\Gw)}\norm{\BBD_\Gw[\xi]}_{L^{q'}(\prt\Gw)}
\leq c\norm{\BBN_\Gw[\gm]}_{L^q(\prt\Gw)}\norm\xi_{W^{1,q'}(\prt\Gw)},
$$
which implies that $\xi\mapsto \langle \gm,\xi\rangle$ is a continuous linear map on $W^{1,q}(\prt\Gw)$, thus belongs to $W^{-1,q}(\prt\Gw)$ and there holds
\bel{IV-1-5}\BA {lll}
\norm{\gm}_{W^{-1,q}(\prt\Gw)}\leq c\norm{\BBN_\Gw[\gm]}_{L^q(\prt\Gw)}.
\EA\ee

\noindent Conversely, if $\gm\in W^{-1,q}(\prt\Gw)$ and $\xi\in C^2(\prt\Gw)$, we set $\phi=\BBN_\Gw[\gm]$ and $w=\BBD_\Gw[\xi]$. Then, using $(\ref{IV-1-4})$
$$\BA {lll}\left|\myint{\prt\Gw}{}\BBN_\Gw[\gm\lfloor_{\prt\Gw}]\BBD_\Gw[\xi]\right|=|\langle \gm,\xi\rangle|\leq \norm\gm_{W^{-1,q}(\prt\Gw)}\norm\xi_{W^{1,q'}(\prt\Gw)}\\\phantom{\left|\myint{\prt\Gw}{}\BBN_\Gw[\gm\lfloor_{\prt\Gw}]\BBD_\Gw[\xi]\right|}
\leq c^{-1}\norm\gm_{W^{-1,q}(\prt\Gw)}\norm{\BBD_\Gw[\xi]}_{L^{q'}(\prt\Gw)}.
\EA$$
By density, it implies that 
 $$\BA {lll}\left|\myint{\prt\Gw}{}\BBN_\Gw[\gm\lfloor_{\prt\Gw}]h\right| \leq c^{-1}\norm\gm_{W^{-1,q}(\prt\Gw)}\norm{h}_{L^{q'}(\prt\Gw)},
 \EA$$
 Hence $\BBN_\Gw[\gm\lfloor_{\prt\Gw}]\in L^q(\prt\Gw)$ and 
 \bel{IV-1-6}\BA {lll}
\norm{\BBN_\Gw[\gm]}_{L^q(\prt\Gw)}\leq c^{-1}\norm{\gm}_{W^{-1,q}(\prt\Gw)}.
\EA\ee
\qeda

\blemma{W} Assume $\gm\in W^{-1,q}(\prt\Gw)\cap\mathfrak M_+(\prt\Gw)$, then problem $(\ref{I-1-23})$ admits a weak solution. 
\es
\Proof We denote by $u_k$ the solution of 
 \bel{IV-1-7}\BA {lll}
\phantom{---}-\Gd u_k+u_k=0\qquad&\text{in }\,\Gw\\
\myfrac{\prt u_k}{\prt{\bf n}}+\min\{u_k^p,k^p\}=\gm&\text{on }\,\prt\Gw,
\EA\ee
the existence of which comes from Theorem D. Then $0\leq  u_k\leq \BBN[\gm]$, and the sequence $\{u_k\}$ is nonincreasing and bounded from above by $\BBN_\Gw[\gm]$. Then it converges to some nonnegative harmonic function $u$ in $\Gw$. Since $\min\{u_k^p,k^p\}\leq (\BBN[\gm])^p$ and 
$\BBN[\gm]\in L^p(\prt\Gw)$, it implies that $\min\{u_k^p,k^p\}$ converges to $u^p$ a.e. in $\prt\Gw$ and in $L^1(\prt\Gw)$. Then $u$ satisfies $(\ref{I-1-23})$.\qeda

\bcor{increa} Let $\{\gm_m\}$ be an increasing sequence of nonnegative measures on $\prt\Gw$ belonging to $W^{-1,p}(\prt\Gw)$ and converging to a measure $\gm$ in 
$\mathfrak M(\Gw)$. Then problem $(\ref{I-1-23})$ with boundary data $\gm$ admits a weak solution. 
\es
\Proof Let $u_m$ be the solution of 
  \bel{IV-1-8}\BA {lll}
-\Gd u_m+u_m=0\qquad&\text{in }\,\Gw\\\phantom{,,,}
\!\myfrac{\prt u_m}{\prt{\bf n}}+u_m^p=\gm_m&\text{on }\,\prt\Gw,
\EA\ee
The sequence $\{u_m\}$ is increasing. For any $\gz\in\CC(\Gw)$, there holds
\bel{IV-1-9}\myint{\Gw}{}u_m\left(-\Gd\gz+\gz\right)dx+\myint{\prt\Gw}{}u^p_m\gz dS=\myint{\prt\Gw}{}\gz d\gm_m
\ee
If we take in particular $\gz=1$, then
$$\myint{\Gw}{}u_mdx+\myint{\prt\Gw}{}u^p_m dS=\myint{\prt\Gw}{}d\gm_m\leq \myint{\prt\Gw}{}d\gm.
$$
Then $u_m$ is bounded in $L^{\frac{N}{N-2},\infty}(\Gw)$ (or any $L^q(\Gw)$ if $N=2$) and $\nabla u_m$ is bounded in $L^{\frac{N}{N-1},\infty}(\Gw)$. By the monotone convergence theorem $\{u_m\}$ converges in $L^1(\Gw)$ to some $u$ and  $\{u_m\lfloor_{\prt\Gw}\}$ converges in $L^p(\prt\Gw)$ to 
$u\lfloor_{\prt\Gw}$. Letting $m\to\infty$ in $(\ref{IV-1-9})$ we obtain
$$\myint{\Gw}{}u\left(-\Gd\gz+\gz\right)dx+\myint{\prt\Gw}{}u^p\gz dS=\myint{\prt\Gw}{}\gz d\gm,
$$
which ends the proof.\qeda\medskip

In the next result we denote by $C^{1,p'}_{_{\BBR^{N-1}}}$ the Bessel (or Sobolev) capacity on $\prt\Gw$ associated to $W^{1,p'}(\BBR^{N-1})$. The 
corresponding capacity  $C^{1,p'}_{{\prt\Gw}}$ on the boundary is defined by local charts and the zero-capacity property does not depend on the charts.  

\bprop {nece} Let $\gm\in\mathfrak M_+(\prt\Gw)$ such that problem $(\ref{I-1-23})$ admits a weak solution, then $\gm$ vanishes on Borel set 
$E\subset\prt\Gw$ satisfying $C^{1,p'}_{\prt\Gw}(E)=0$. 
\es
\Proof Without loss of generality we can assume that $E\subset\prt\Gw$ is a compact set. Because of uniqueness, $u$ is nonnegative. Let $\eta\in C_0^2(\prt\Gw)$ such that $0\leq\eta\leq 1$, $\eta=1$ in a neighborhood of $E$ and $v_\eta=\BBP_{\Gw}[\eta]$. If $\gz\in C^2(\overline\Gw)$ we have
$$\myint{\Gw}{}u\left(-\Gd\gz+\gz\right)dx+\myint{\prt\Gw}{}u^p\gz dS=-\myint{\prt\Gw}{}\myfrac{\prt\gz}{\prt{\bf n}} udS+\myint{\prt\Gw}{}\gz d\gm.
$$
For $k>1$ we take $\gz=v_\eta^k$, then 
$$\BA{lll}\myint{\Gw}{}u\left(-k(k-1)v_\eta^{k-2}|\nabla v_\eta|^2+(1-k)v_\eta^k\right)dx
+\myint{\prt\Gw}{}u^pv_\eta^k dS\\[2mm]\phantom{--------------}
=-k\myint{\prt\Gw}{}\eta^{k-1} \BBD_\Gw[\eta]udS +\myint{\prt\Gw}{}\eta^k d\gm.
\EA$$
Since $v_\eta\geq 0$, we obtain
$$\myint{\prt\Gw}{}u^pv_\eta^k dS+k\myint{\prt\Gw}{}\eta^{k-1} |\BBD_\Gw[\eta]|udS\geq\gm(E).
$$
Furthermore.
$$\myint{\prt\Gw}{}\eta^{k-1} |\BBD_\Gw[\eta]|udS\leq \left(\myint{\prt\Gw}{}\eta^{k} u^pdS\right)^{\frac1p}
\left(\myint{\prt\Gw}{}\eta^{k-p'} |\BBD_\Gw[\eta]|^{p'}dS\right)^{\frac1{p'}}.
$$
Taking $k=p'$ and using $(\ref{IV-1-4})$, we infer
\bel{IV-1-10}\myint{\prt\Gw}{}u^pv_\eta^k dS+cp'\left(\myint{\prt\Gw}{}\eta^{k} u^pdS\right)^{\frac1p}\norm\eta_{W^{1,p'}(\prt\Gw)}\geq\gm(E).
\ee
If $C^{1,p'}_{\prt\Gw}(E)=0$, there exists a sequence $\{\eta_m\}\subset C_0^2(\prt\Gw)$ such that $0\leq\eta_m\leq 1$, $\eta=1$ in a neighborhood of $E$ and 
$\norm{\eta_m}_{W^{1,p'}(\prt\Gw)}\to 0$ when $m\to\infty$. This implies that $v_{\eta_m}\to 0$ in $L^1(\Gw)$, hence the left-hand side of 
$(\ref{IV-1-10})$ tends to $0$, and finally $\gm(E)=0$.\qeda\medskip

We end the proof of Theorem F with the sufficient condition which follows from a general result due to Feyel and de la Pradelle \cite{FeyPra}. 
\bprop {suf} Let $\gm\in\mathfrak M_+(\prt\Gw)$ such that $\gm(E)=0$ for any Borel set 
$E\subset\prt\Gw$ satisfying $C^{1,p'}_{\prt\Gw}(E)=0$. Then there exists an increasing sequence $\{\gm_n\}\subset\mathfrak M_+(\prt\Gw)\cap W^{-1,p}(\partial\Gw)$ converging to $\gm$.  
\es

\noindent \Remark If $1<p<\frac{N-1}{N-2}$, $W^{1,p'}(\prt\Gw)\subset C(\prt\Gw)$. Therefore the only set with zero $C^{1,p'}_{\prt\Gw}$-capacity is the empty set. If $p\geq \frac{N-1}{N-2}$, a single point has zero $C^{1,p'}_{\prt\Gw}$-capacity. Since $\gd_a(a)=1$ for any $a\in\prt\Gw$ there is no solution of problem $(\ref{I-1-23})$ with $\gm=\gd_a$. \medskip

As a consequence we have a non-removability result.

\bcor{removK} Let $\Gw$ be a smooth bounded domain of $\BBR^N$. Then any compact subset $K\subset\prt\Gw$ with positive $C^{1,p'}_{{\prt\Gw}}$-capacity is non-removable in the sense that there exists a nonnegative non-trivial function $u_K\in C^{1}(\overline\Gw\setminus K)$ satisfying 
  \bel{IV-1-11}\BA {lll} 
\phantom{|u^{p-1}u,}
-\Gd u+ u=0\qquad&\text{in }\;\Gw\\[0mm]\phantom{+ u,,}
\myfrac{\prt u}{\prt {\bf n}}+ |u|^{p-1}u=0\qquad&\text{in }\;\prt\Gw\setminus K.
\EA
\ee
\es
\Proof By \cite[Theorem 2.5.3]{AH} there exists a positive measure, called the capacitary measure $\gm_K$ with support in $K$ and such that 
$\gm_K\in W^{-1,p'}(\Gw)$. For such a measure there exists a positive solution to $(\ref{I-1-23})$, hence $u$ satisfies $(\ref{IV-1-11})$. 
\phantom{zh"(qkjÜʲ×}\qeda
\medskip

 \noindent \Remark We conjecture that the condition $C^{1,p'}_{\prt\Gw}(K)=0$ is also a sufficient condition for a compact set $K\subset\prt\Gw$ to be removable. This is even not known if $K$ is a singleton.

\bibliographystyle{plain}

\end{document}